\documentclass[11pt,twoside]{article}

\usepackage{fullpage}

\usepackage{epsf}
\usepackage{fancyhdr}
\usepackage{graphics}
\usepackage{graphicx}
\usepackage{psfrag}
\usepackage{algorithmic}

\usepackage[linesnumbered,ruled]{algorithm2e}
\DontPrintSemicolon

\usepackage{color}

\usepackage{amsthm}
\usepackage{amsfonts}
\usepackage{amsmath}
\usepackage{amssymb,bbm}
\usepackage{natbib}
\usepackage{nicefrac}
\usepackage{url} 
\usepackage[colorlinks=True,linkcolor=magenta,citecolor=blue,urlcolor=blue,pagebackref=true,backref=true]
{hyperref}
\renewcommand*{\backref}[1]{\ifx#1\relax \else Page #1 \fi}
\renewcommand*{\backrefalt}[4]{%
    \ifcase #1 \footnotesize{(Not cited.)}%
    \or        \footnotesize{(Cited on page~#2.)}%
    \else      \footnotesize{(Cited on pages~#2.)}%
    \fi}

\usepackage{chngpage}
\usepackage{tabularx}
\usepackage{enumitem}
\usepackage{booktabs}
\usepackage{caption,subcaption}
\usepackage{mathtools}
\usepackage{afterpage}
 
\usepackage{tikz}
\usetikzlibrary{calc, shapes, backgrounds}
\usepackage{subcaption} 

\numberwithin{equation}{section}
\theoremstyle{plain}

\newtheorem{theorem}{Theorem}[section]
\newtheorem{lemma}[theorem]{Lemma}

\newtheorem{proposition}[theorem]{Proposition}

\usepackage{color}
\usepackage{amssymb}
\usepackage{comment}
\usepackage{pdfpages}
\usepackage{graphicx}
\usepackage{dcolumn}
 
\numberwithin{equation}{section}

\newcommand{\parenth}[1]{\left( #1 \right)}

\newcommand{\abss}[1]{\left| #1 \right |}

\newcommand{\ordersym}{\bar{r}_{\text{sym}}}
\newcommand{\orderassym}{\bar{r}_{\text{asym}}}

\newtheoremstyle{named}{}{}{\itshape}{}{\bfseries}{.}{.5em}{\thmnote{#3's }#1}
\theoremstyle{named}

\theoremstyle{plain}

\newlength{\widebarargwidth}
\newlength{\widebarargheight}
\newlength{\widebarargdepth}

\makeatletter
\long\def\@makecaption#1#2{
        \vskip 0.8ex
        \setbox\@tempboxa\hbox{\small {\bf #1:} #2}
        \parindent 1.5em  
        \dimen0=\hsize
        \advance\dimen0 by -3em
        \ifdim \wd\@tempboxa >\dimen0
                \hbox to \hsize{
                        \parindent 0em
                        \hfil
                        \parbox{\dimen0}{\def\baselinestretch{0.96}\small
                                {\bf #1.} #2
                                }
                        \hfil}
        \else \hbox to \hsize{\hfil \box\@tempboxa \hfil}
        \fi
        }
\makeatother

\long\def\comment#1{}
\definecolor{battleshipgrey}{rgb}{0.52, 0.52, 0.51}
\definecolor{darkgray}{rgb}{0.66, 0.66, 0.66}
\definecolor{darkgreen}{rgb}{0.0, 0.2, 0.13}
\definecolor{darkspringgreen}{rgb}{0.09, 0.45, 0.27}
\definecolor{dukeblue}{rgb}{0.0, 0.0, 0.61}
\definecolor{olivedrab7}{rgb}{0.24, 0.2, 0.12}
\definecolor{darkblue}{rgb}{0.0, 0.0, 0.55}
\definecolor{darkscarlet}{rgb}{0.34, 0.01, 0.1}
\definecolor{candyapplered}{rgb}{1.0, 0.03, 0.0}
\definecolor{ao(english)}{rgb}{0.0, 0.5, 0.0}
\definecolor{applegreen}{rgb}{0.55, 0.71, 0.0}

\usepackage{macros}

\begin{document}

\begin{center}
\textbf{\Large{Uniform Convergence Rates for Maximum Likelihood \\[0.2pt] Estimation under Two-Component Gaussian Mixture Models}}
\end{center}
 
 {\large{
 \begin{center}
 \begin{tabular}{cc}
 Tudor Manole$^{\star, \diamond}$ & Nhat Ho$^{\star, \dagger}$ 

 \end{tabular}
 \end{center}
}}

 {\large{
 \begin{center}
 \begin{tabular}{cc}
Department of Statistics and Data Science$^\diamond$ \\
Carnegie Mellon University\\[0.1in] 
Department of Electrical Engineering and Computer Science$^\dagger$\\
University of California, Berkeley 
 \end{tabular}
 \end{center}
}}


\begin{abstract}
We derive uniform convergence rates for the maximum likelihood estimator and minimax lower bounds for parameter
estimation in two-component
location-scale Gaussian mixture models with unequal variances.
We assume the mixing proportions of the mixture are known and fixed, but make no separation
assumption on the underlying mixture components. 
A phase transition is shown to exist in the optimal parameter estimation rate, depending
on whether or not the mixture is balanced. Key to our analysis is a careful study of the 
dependence between the parameters of location-scale
Gaussian mixture models, as captured through systems of polynomial equalities and inequalities
whose solution set drives the rates we obtain. 
A simulation study illustrates the theoretical findings of this work.
\end{abstract}
\let\thefootnote\relax\footnotetext{$\star$ Tudor Manole and Nhat Ho contributed equally to this work.}

\section{Introduction}
\label{sec:intro}
Finite mixture models are a widely-used tool for modeling heterogeneous data,
consisting of hidden  subpopulations
with distinct distributions. For applications exhibiting 
continuous data, 
location-scale Gaussian mixtures are arguably the most popular family of parametric
mixture models. 
Beyond their broad applications as a modeling and clustering tool in the social, physical and life sciences
\citep{mclachlan2004}, Gaussian mixtures provide
a flexible approach to density estimation \citep{genovese2000,ghosal2001}.

Estimating the parameters of a mixture model is crucial
for quantifying the underlying heterogeneity of the data. One of the most widely-used approaches is the maximum likelihood
estimator (MLE). 
A Gaussian mixture model with a known number of components $K$, all of which are well-separated,
forms a regular parametric model for which the MLE achieves the standard parametric estimation rate
\citep{ho2016c,chen2017a}. 
Such  rates are typically understood in terms of convergence of mixing measures, quantified
using the Wasserstein distance as a means of avoiding label switching issues
inherent in mixture modeling \citep{nguyen2013a}.
In the absence of separation conditions, mixture components are permitted to overlap
arbitrarily, thus the number of distinct components, say $K_0$, may be strictly less than $K$. In this setting,
the Fisher information matrix of the mixture model becomes singular, and has been shown to lead to slower rates
of paramater estimation. For instance, \cite{ho2016b} showed that the pointwise convergence rate
of the MLE under location-scale Gaussian mixtures deteriorates as the difference $K-K_0$ increases. 
Here, the term ``pointwise'' refers to the rates therein being dependent upon the parameters of the true underlying
mixture. These rates therefore do not provide upper bounds on the worst-case risk, and hence on the minimax risk.
To the best of our knowledge, minimax rates for parameter estimation in general location-scale Gaussian mixtures have only
been studied by \cite{hardt2015} in the case $K=2$, using estimators different than the MLE.

\textbf{Our Contributions.} In this paper, we establish 
uniform convergence rates of the MLE under one-dimensional, two-component location-scale Gaussian mixture models with unequal variances.
Our rates differ substantially from the pointwise rates of \cite{ho2016b}. 
We show that the optimal estimation rate differs 
according to whether or not the underlying mixture has equal mixing proportions, which
we refer to as a symmetric mixture.
This phase transition motivates us to restrict our analysis to mixtures admitting fixed and known
mixing proportions. We also prove that these rates are minimax optimal, thereby
refining the known minimax rates from \cite{hardt2015} to the distinct settings of symmetric and asymmetric two-component
mixtures. 
Our analysis relies upon the strong dependence between the parameters of location-scale
Gaussian mixture models. Indeed, the rates we obtain are driven by the solution set
of explicit systems of polynomial equalities and inequalities, arising from a key
linear dependence between certain partial derivatives of Gaussian densities, described in equation \eqref{eq:key_pde} below.

\subsection{Related literature\label{sec:related}}
Establishing optimal rates for parameter estimation in finite mixture models
is a long-standing problem, dating back at least to the seminal work of \cite{chen1995}. 
For one-dimensional mixtures with a number of components $K_0$ which is unknown but bounded
above by a known constant $K$, 
\cite{chen1995} showed that the optimal pointwise rate of parameter estimation 
scales as $C_0 n^{-1/4}$, where $n$ is the sample size, and $C_0 > 0$ 
is a constant depending on the underlying true mixture model in a possibly unbounded manner. This result holds
for mixtures satisfying a condition known as strong identifiability, which requires the mixture component densities
and a certain number of their partial derivatives to be linearly independent---a condition satisfied
by location Gaussian mixtures, but not location-scale Gaussian mixtures.
\cite{nguyen2013a} and \cite{ho2016c} also establish the $C_0 n^{-1/4}$ pointwise rate for multivariate
strongly identifiable mixtures with fixed dimension. These pointwise rates do not, however, provide upper
bounds on the minimax risk of parameter estimation, due to the lack of uniformity in the constant $C_0$. 
Indeed, \cite{heinrich2018} proved that this minimax risk, under strongly identifiable mixtures,
scales at the markedly distinct rate
$n^{- \frac 1 {4(K-K_0) + 2}}$, which deteriorates exponentially with the level of overspecification $K-K_0$
of the number of components. In this context, the quantity $K_0$ is understood as the minimum number of well-separated 
components
of the underlying mixture, with the case $K_0=1$ corresponding to the rate with no separation assumption.
The minimax rate established by \cite{heinrich2018} is achievable by a minimum-distance estimator, 
and by the Denoised Method of Moments \citep{wu2019}. A multivariate extension of the latter method was also shown
to achieve the minimax rate of estimating a high-dimensional location-Gaussian mixture model \citep{doss2020}---see also \cite{wu2019a} for the special case $K=2$ of the minimax rate therein. We refer to 
\cite{vempala2004,moitra2010, kalai2010a, azizyan2013} and references therein for prior advances in the high-dimensional setting. 

For mixture models failing to satisfy the strong identifiability condition, 
optimal rates for parameter estimation do not enjoy a unified treatment. For Gaussian mixture models with unknown means and common but unknown variances, 
\cite{wu2019} showed that the $n^{-\frac 1 {4K}}$ rate is minimax optimal under no separation assumptions, and achievable
by the Denoised Method of Moments. \cite{feller2019} shows this rate is also achievable by the MLE under two-component mixtures with equal variances. 
When the variances of the Gaussian mixture are allowed to be unknown and distinct, \cite{ho2016b} establish the
pointwise rate $C_0 n^{-\frac 1 {2r}}$, for an integer $r \geq 1$ determined by the solution set of a system of polynomial equations depending on $K$---for
instance, one has $r=2$ in the case $K=2$, leading to the $C_0 n^{-1/8}$ pointwise rate which had previously been observed by \cite{chen2003}. 
In contrast, the minimax rate in the two-component case 
was shown to be $n^{-1/12}$ by \cite{hardt2015} using the method of moments. As we will show in this paper, 
the minimax rate $n^{-1/12}$ can be improved to $n^{-1/8}$ for symmetric mixtures, up to a polylogarithmic factor. 
To the best of our knowledge, 
our work is the first to provide uniform upper bounds on the rate of convergence of the MLE, and we show that 
it is minimax optimal both for asymmetric and symmetric mixtures respectively, up to polylogarithmic factors.



\subsection{Problem Setting}
\textbf{Gaussian Mixture Models and Maximum Likelihood Estimation.} 
Throughout this paper, we fix two compact subsets  $\Theta$ and $\Omega$ of $\mathbb{R}$ and $\mathbb{R}_{+}$ respectively,
such that $0 \in \text{int}(\Theta)$, where $\text{int}(\cdot)$ denotes the interior of a set.  
Let $\calF = \{f(\cdot,\theta, \sigma^2): \theta \in \Theta, \sigma^2 \in \Omega\}$ denote the 
location-scale Gaussian parametric family, where 
$$f(x, \theta, \sigma^2) = \frac 1 {\sqrt{2\pi\sigma^2}} \exp\left(-\frac{(x-\theta)^2}{2\sigma^2}\right),
\quad x \in \bbR.$$
Fix a known real number $\pi \in (0,1/2]$, and let $c = \pi/(1-\pi)$. 
Let $Y_1, \dots, Y_n$ be an i.i.d. sample from the the one-dimensional
location-scale Gaussian mixture model whose density is given by
\begin{align}
\label{eqn:general_model}
g(x, \bfeta_n) = \pi f(x,-\theta_{n},\sigma_{1,n}^2) + (1-\pi)f(x,c\theta_{n},\sigma_{2,n}^2),\quad x \in \bbR
\end{align}
where $\bfeta_n=(\theta_n, \sigma_{1,n}^2, \sigma_{2,n}^2) \in \etaspace$, 
and $\etaspace = \Theta \times \Omega^2$. 
We will also use the shorthand $v_{n,j} = \sigma^2_{n,j}$ for $j=1,2$, in the sequel.
We focus on model \eqref{eqn:general_model} throughout the paper.
To emphasize the uniformity in our bounds below, notice that we allow for the parameters 
$\bfeta_n$ to vary with the sample size $n$, converging to some limit points.
Notice further that the choice of parametrization in model \eqref{eqn:general_model} 
ensures that the mixture model has zero mean. Our results can be extended to mixtures with general mean
$\mu \in \bbR$, whose density is of the form
$\pi f(\cdot,\mu-\theta_{n},\sigma_{1,n}^2) + (1-\pi)f(\cdot,\mu+c\theta_{n},\sigma_{2,n}^2)$,
but we only consider the case $\mu=0$ for simplicity.

The log-likelihood function of $\bfeta_n$ with respect to the sample $Y_1, \dots, Y_n$ is given by
$$\ell_n(\bfeta) = \sum_{i=1}^n \log g(Y_i; \bfeta), \quad \bfeta \in H.$$
We let $\hbfeta_n = (\htheta_n, \hsigma_{1,n}^2, \hsigma_{2,n}^2)$ denote a maximizer of $\ell_n$
over $\etaspace$. The existence of $\hbfeta_n$ is guaranteed by the compactness of the parameter space $H$.

\textbf{Loss Function on $H$.} In order to quantify the convergence of parameters in $H$, we introduce the following loss function.
Given $\bfeta^{(1)} = (\theta^{(1)}, v_1^{(1)}, v_2^{(1)}), \bfeta^{(2)} = (\theta^{(2)}, v_1^{(2)}, v_2^{(2)}) \in \etaspace$,
define 
\begin{align}
\symloss_r (\bfeta^{(1)},\bfeta^{(2)})
 : = \min \biggr\{&|\theta^{(1)} - \theta^{(2)}|^r+|v_{1}^{(1)} - v_{1}^{(2)}|^{r/2} + |v_{2}^{(1)} - v_{2}^{(2)}|^{r/2}, \nonumber \\
&                |\theta^{(1)} + \theta^{(2)}|^r+|v_{1}^{(1)} - v_{2}^{(2)}|^{r/2} + |v_{2}^{(1)} - v_{1}^{(2)}|^{r/2}\biggr\}^{\frac 1 r}, 
\end{align}
where $r \geq 1$. Notice that $\symloss_r$ is invariant to label switching of mixture components, 
and reduces to the loss function used by \cite{hardt2015} in the special case $r=1$.
To understand how convergence under $\symloss_r$ relates to convergence of the individual
mixture parameters, let $\bar\bfeta_n = (\bar\theta_n, \bar v_{1,n}, \bar v_{2,n}) \in H, n \geq 1,$ be a sequence satisfying
$\symloss_r(\bar\bfeta_n,\bfeta_n) \leq \alpha_n$, for a sequence 
of nonnegative real numbers $\alpha_n \to 0$. Then there exists a permutation $\tau$ on $\{1,2\}$ such that
$$\Big| |\bar \theta_n^{(1)}| - |\theta_n^{(2)}| \Big| \lesssim \alpha_n, \quad
  |\bar v_{j,n}^{(1)} -  v_{\tau(j),n}^{(2)}| \lesssim \alpha_n^2, \quad j=1,2.$$
The loss function $\varphi_r$ captures, in particular, the inhomogeneity in estimating
the means and variances of a Gaussian mixture model---indeed, it has been observed at least since
the work of \cite{chen2003} that typical rates of convergence for the variances of a Gaussian mixture are
faster than those of their means. 
We also note that $\varphi_r$ admits a natural interpretation in terms of the Wasserstein distance, 
a metric frequently used for quantifying convergence rates in multivariate 
mixtures with more than two components \citep{nguyen2013a, heinrich2018}. 
Specifically, defining probability measures $\bar G_n = \pi \delta_{-\bar \theta_n} + 
(1-\pi) \delta_{c\bar \theta_n}$, $G_n =  \pi \delta_{- \theta_n} + 
(1-\pi) \delta_{c\theta_n}$,
and $\bar H_n = \pi\delta_{\bar v_{1,n}^{(1)}} + (1-\pi) \delta_{\bar v_{2,n}^{(j)}}$,
$H_n = \pi\delta_{v_{1,n}} + (1-\pi) \delta_{v_{2,n}}$, where $\delta_x$ denotes a Dirac
measure placing mass at $x \in \bbR$, we have
$$\symloss_r^r(\bar\bfeta_n, \bfeta_n) \asymp W_r^r(\bar G_n, G_n) + 
											  W_{r/2}^{r/2}(\bar H_n, H_n),$$
where $W_r$ denotes the $r$-th order Wasserstein distance (see~\cite{villani2003} for a formal definition of the Wasserstein distance).

Finally, since we have assumed in model \eqref{eqn:general_model} that the mixing proportion
$\pi$ is known and fixed, mixture label switching generically occurs only in the 
symmetric setting $\pi=1/2$. 
When working in the asymmetric setting $\pi\neq 1/2$ below, we will therefore be able to state our
results in terms of the stronger loss function
$$\assymloss_r(\bfeta^{(1)}, \bfeta^{(2)}) = \Big(|\theta^{(1)} - \theta^{(2)}|^r+|v_{1}^{(1)} - v_{1}^{(2)}|^{r/2} + |v_{2}^{(1)} - v_{2}^{(2)}|^{r/2}\Big)^{\frac 1 r},$$
for all $r \geq 1$.

\subsection{Paper Outline}
The rest of this paper is organized as follows. 
In Section \ref{sec:main}, we state our main results regarding the rate
of convergence of the MLE and minimax lower bounds, both in the asymmetric regime (Section \ref{sec:asymmetric})
and the symmetric regime (Section \ref{sec:symmetric}). In Section \ref{sec:simulations}
we illustrate our theoretical findings with a simulation study. We close with discussions in Section 
\ref{sec:discussion}. All proofs are relegated to Appendices \ref{sec:appendix_asymmetric},
\ref{sec:appendix_symmetric} and \ref{sec:appendix_rasym}, 
and further simulation specifications are included in Appendix
\ref{sec:appendix_numerical}. 

\subsection{Notation}
For any two densities $p$ and $q$ with respect to Lebesgue measure, the Total Variation
distance between $p$ and $q$ is given by $V(p,q) = (1/2)\int \abss{p-q}d\mu$, and 
the squared Hellinger distance between $p$ and $q$ is 
given by $h^{2}(p,q)= (1/2)\int \parenth{p^{1/2}-q^{1/2}}^2d\mu$. 
Given two sequences of nonnegative real numbers $(a_n)_{n=1}^\infty, (b_n)_{n=1}^\infty$,
we write $a_{n} \gtrsim b_{n}$ if there exists a constant $C > 0$
not depending on $n$
such that $a_{n} \geq C b_{n}$ for all $n \geq 1$. We write $a_n \asymp b_n$ if $a_n \lesssim b_n \lesssim a_n$.
For any multi-index $\alpha=(\alpha_1, \dots, \alpha_k)$ where $\alpha_1,\dots,\alpha_k \in \bbN$,
we write $|\alpha| = \sum_{i=1}^k \alpha_i$.

\section{Convergence Rates of the Maximum Likelihood Estimator and Minimax Lower Bounds}
\label{sec:main}
In this section, we state our main results regarding the uniform rate of convergence of the MLE
and corresponding minimax lower bounds. 
Key to our analysis is a careful treatment of the dependence 
between the mean and variance parameters of model \eqref{eqn:general_model}, 
which is determined by the following partial differential equation (PDE)
satisfied by the Gaussian density $f$,
\begin{eqnarray}
\label{eq:key_pde}
\frac{\partial^{2}{f}}{\partial{\theta^{2}}}(x,\theta,v) = 2 \frac{\partial{f}}{\partial{v}}(x,\theta,v), \quad
x \in \bbR, ~ \theta \in \Theta,~  v \in \Omega.
\end{eqnarray}
This equality prevents location-scale Gaussian densities from satisfying the strong
identifiability criterion, for which minimax rates are well understood \citep{heinrich2018}, 
and will lead to worse rates of convergence for parameter estimation in the sequel. 
Under the specific setting that we consider, equation \eqref{eq:key_pde} 
also creates a new phase transition in the parameter estimation rates,
under the two regimes $\pi \in (0,1/2)$ and $\pi = 1/2$, which have not been addressed so far in the literature.
We treat these two regimes separately in what follows.

\subsection{Asymmetric Regime}
\label{sec:asymmetric}
 
Throughout this subsection, we assume $\pi \in (0,1/2)$ is known and fixed. 
The convergence rate of the MLE under the asymmetric regime is governed by the solution
set of a system of polynomial equations which we now describe. 
Let $\orderassym(\pi)$ denote 
 the smallest positive integer $r \geq 1$ such that the following system of polynomial equations
\begin{align}
\hspace{ - 3 em} (1-\pi) \sum \limits_{\alpha_{1},\alpha_{2},\beta_{1},\beta_{2}} \dfrac{1}{2^{\alpha_{2}+\beta_{2}}} \dfrac{c^{\alpha_{1}}(c+1)^{\beta_{1}}(-x_{1})^{\alpha_{1}}x_{2}^{\beta_{1}}y_{2}^{\alpha_{2}}y_{3}^{\beta_{2}}}{\alpha_{1}!\alpha_{2}!\beta_{1}!\beta_{2}!} \nonumber \\
+ \pi \sum \limits_{\alpha_{1},\alpha_{2}} \dfrac{1}{2^{\alpha_{2}}}\dfrac{x_{1}^{\alpha_{1}}y_{1}^{\alpha_{2}}}{\alpha_{1}!\alpha_{2}!}  = 0, ~~
 \text{for each}\ \ell=1,\ldots,r \label{eq:asymmetric_system}
\end{align}
does not have any non-trivial real-valued solution for $(x_{1},x_{2},y_{1},y_{2},y_{3}) \in \bbR^5$. Here, the range of 
the first sum is over all nonnegative integers $\alpha_{1},\alpha_{2},\beta_{1},\beta_{2}$ such that
$\alpha_{1}+\beta_{1}+2\alpha_{2}+2\beta_{2} = \ell$, $1 \leq \alpha_{1}+\alpha_{2} \leq r$, 
and $0 \leq \beta_{1}+\beta_{2} \leq r - (\alpha_{1}+\alpha_{2})$,
while the ranges of $\alpha_{1},\alpha_{2}$ in the second sum satisfy 
$\alpha_{1}+2\alpha_{2}=\ell$, $1 \leq \alpha_{1}+\alpha_{2} \leq r$. 
A solution is considered non-trivial if at least one of the variables $x_{1}, y_{1}$, and $y_{2}$ is different from 0.
The quantity $\orderassym(\pi)$ is called the asymmetric order, 
and we now show its central role in the convergence rate of the MLE under the asymmetric
regime. 
\begin{theorem}  \label{theorem:convergence_asym}
Let $\pi \in (0,1/2)$ be fixed.
\begin{itemize}
\item[(a)] (Maximum Likelihood Estimation) We have
\begin{align*}
 \sup_{\bfeta \in H} 
\bbE_{\bfeta}\Big[\assymloss_{\orderassym(\pi)}(\hbfeta_n,\bfeta)\Big] \lesssim \left(\frac{\log n}{n}\right)^{\frac 1 {2\orderassym(\pi)}}
 \end{align*}
where the expectation is
taken with respect to the product distribution of an i.i.d. sample $Y_{1},\ldots,Y_{n}$ 
from model \eqref{eqn:general_model}.
\item[(b)] (Minimax Lower Bound) Let $v_0 \in \text{int}(\Omega)$,  $\bfeta_0 = (0,v_0, v_0)$, and
$$\etaspace(\kappa) = \left\{\bfeta \in \etaspace: \assymloss_{\orderassym(\pi)}^{\orderassym(\pi)}(\bfeta,\bfeta_0) \leq \kappa \right\}, \quad \kappa > 0.$$
 Then,
there exists a universal constant $c_1 > 0$ such that
\begin{align*}
\inf_{\tilde\bfeta_n}
\sup_{\substack{\bfeta \in \etaspace(c_1n^{-1/2})}} 
\bbE_{\bfeta}\Big[\assymloss_{\orderassym(\pi)}(\tilde\bfeta_n,\bfeta)\Big] \gtrsim \left(\frac{1}{n}\right)^{\frac 1 {2\orderassym(\pi)}},
\end{align*}
where the infimum is over all sequences of estimators $\tilde \bfeta_n$ based on an i.i.d. sample $Y_1, \dots, Y_n$ from 
model \eqref{eqn:general_model}.
\end{itemize}
\end{theorem} 
Theorem \ref{theorem:convergence_asym}(a)
implies that the rate of convergence of the MLE
under $\assymloss_{\orderassym(\pi)}$ is of order
$n^{-1/2\orderassym(\pi)}$ up to a polylogarithmic factor. 
To prove this result, our key theoretical contribution is a characterization of
the distance $\assymloss_r$ between mixture parameters
in terms of the Total Variation distance between their corresponding mixture densities, 
a general approach which has previously formed the basis minimax analyses
for strongly identifiable mixture models \citep{heinrich2018, doss2020}
and pointwise parameter estimation rates
for location-scale Gaussian mixture models \citep{ho2016b}. 
Specifically, we prove in Theorem \ref{theorem:prelim_asym} in Appendix \ref{sec:appendix_asymmetric} that
for any $\bfeta^{(1)}, \bfeta^{(2)} \in H$, 
the following inequality holds
\begin{equation}
\label{eq:psi_tv}
\assymloss_{\orderassym(\pi)}(\bfeta^{(1)}, \bfeta^{(2)}) \lesssim V\Big(g(\cdot,\bfeta^{(1)}), g(\cdot, \bfeta^{(2)})\Big)^{\frac 1 {\orderassym(\pi)}}.
\end{equation}
Combining this bound with a generic convergence result for the maximum likelihood
density estimator \citep{vandegeer2000}, together with bracketing entropy bounds for classes of mixture densities
\citep{ghosal2001}, readily
leads to Theorem \ref{theorem:convergence_asym}(a).  Theorem \ref{theorem:convergence_asym}(b)
further shows that the resulting rate is minimax optimal. We wish to emphasize
that the bound \eqref{eq:psi_tv} may similarly be used to obtain convergence rates
for parameter estimation of any other method admitting a known density estimation guarantee.

In order to obtain a quantitative rate of
convergence, we bound the asymmetric order as follows.
\begin{proposition}
\label{proposition:asymmetric_system}
Under the asymmetric system \eqref{eq:asymmetric_system}
with any $\pi \in (0,1/2)$, we have $\orderassym(\pi) \geq  6$.
\end{proposition}
Proposition 
\ref{proposition:asymmetric_system} provides a lower bound on the asymmetric order, which
is obtained through an explicit solution to the asymmetric system \eqref{eq:asymmetric_system}
when $r=5$. Upper bounding $\orderassym(\pi)$ requires showing that the asymmetric system admits no non-trivial solutions
for a given $r \geq 6$, a problem which may be solved using various techniques from algebraic geometry, 
such as the method of Gr\"obner bases  \citep{buchberger1985, sturmfels2005}. 
In Appendix \ref{sec:appendix_rasym}, we apply this method 
with the Mathematica programming language \citep{wolfram1999}  
to show that for all $\pi\in\{i/100:1 \leq i \leq 49, i \in \bbN\}$, 
the system of polynomials \eqref{eq:asymmetric_system}
with $r=6$ admits no solutions. For these values of $\pi$, it follows
that  $\orderassym(\pi) = 6$, and we conjecture this result to hold uniformly over all $\pi \in (0,1/2)$,
but we do not have a proof. 
For the values of $\pi$ where this result holds, 
Theorem \ref{theorem:convergence_asym}(a)
leads to the 
following rates for parameter estimation under model \eqref{eqn:general_model}
\begin{equation}
\label{eq:rate_asymmetric_rasym}
|\htheta_n - \theta_n| \lesssim \left(\frac{\log n}{n}\right)^{1/12}, \quad 
  |\hat v_{j,n} - v_{j,n}| \lesssim \left(\frac{\log n}{n}\right)^{1/6}, ~~ j=1,2,
\end{equation}
with probability tending to one, as $n \to \infty$. 
Equation \eqref{eq:rate_asymmetric_rasym} exhibits a discrepancy between
the convergence rates of the location and scale parameters of the mixture, 
which essentially arises from the key PDE \eqref{eq:key_pde}.

Under the regime where the variances $v_{1,n} = v_{2,n}$ are assumed to be equal
but unknown and $\pi \neq 1/2$, \cite{feller2019} previously
established the $n^{-1/6}$ uniform rate of convergence of the MLE
under the $\assymloss_3$ loss function, up to polylogarithmic factors. Our results imply the significantly slower rate
$n^{-1/12}$, when the variances are not constrained to be equal. Under no assumptions on $\pi$, this
rate was already known to be minimax optimal from \cite{hardt2015}, 
who prove the lower bound
$$\inf_{\tilde \bfeta_n} \sup_{\bfeta \in H} \bbE_{\bfeta} \big[\assymloss_1(\tilde \bfeta_n, \bfeta)\big] \gtrsim n^{-\frac 1 {12}}$$
Noting that $\assymloss_1 \gtrsim \assymloss_r^r$ for all $r \geq 1$, 
our Theorem \ref{theorem:convergence_asym}
recovers this lower whenever $\orderassym(\pi)=6$, using a distinct proof technique, and shows
that it is achievable by the MLE. In 
constrast to the lower bound of \cite{hardt2015}, however, 
we show in what follows that parameter estimation rates are markedly different
in the symmetric regime $\pi = 1/2$.

\subsection{Symmetric Regime}
 \label{sec:symmetric}
We now establish parameter estimation rates in the symmetric regime
where $\pi=1/2$. Unlike the previous subsection, our results will now
be driven by the solution set to a system of both polynomial equations and a polynomial inequality.
Specifically, we denote 
by $\ordersym$ the smallest positive integer $r$ such that the following system 
\begin{align}
& \sum \limits_{\alpha_{1},\alpha_{2},\beta_{1},\beta_{2}} \dfrac{1}{2^{\alpha_{2}+\beta_{2}}} \dfrac{2^{\beta_{1}}(-x_{1})^{\alpha_{1}}x_{2}^{\beta_{1}}y_{2}^{\alpha_{2}}y_{3}^{\beta_{2}}}{\alpha_{1}!\alpha_{2}!\beta_{1}!\beta_{2}!} \nonumber \\
& \hspace{10 em}  + \sum \limits_{\alpha_{1},\alpha_{2}} \dfrac{1}{2^{\alpha_{2}}}\dfrac{x_{1}^{\alpha_{1}}y_{1}^{\alpha_{2}}}{\alpha_{1}!\alpha_{2}!}  = 0 \ \text{for each} \ l=1,\ldots,r,
\label{eq:symmetric_system_equalities} \\
& \abss{x_{1}}^{r}+\abss{y_{1}}^{r/2}+\abss{y_{2}}^{r/2} \leq \abss{2x_{2}-x_{1}}^{r}+\abss{y_{3}-y_{1}}^{r/2}+\abss{y_{2}+y_{3}}^{r/2} \label{eq:symmetric_system_inequalities}
\end{align}
does not admit any non-trivial real-valued solution for $(x_{1},x_{2},y_{1},y_{2},y_{3})\in\bbR^5$. 
Here, the ranges of $\alpha_{1},\alpha_{2},\beta_{1},\beta_{2}$ in the 
above sums as well as the notion of non-triviality are defined similarly as those of
the asymmetric system of polynomial equations \eqref{eq:asymmetric_system}. 
$\ordersym$ is called the symmetric order.
%

We note that the polynomial equalities of the asymmetric system \eqref{eq:asymmetric_system}
reduce to those of the symmetric system \eqref{eq:symmetric_system_equalities} when $\pi=1/2$. On the other hand, 
the present system also contains the inequality \eqref{eq:symmetric_system_inequalities},
which arises due to
the symmetric structure of model \eqref{eqn:general_model} when $\pi=1/2$. The following
straightforward result suggests the necessity of this inequality. 
\begin{proposition}
\label{proposition:inequality_necessity}
For every integer $r \geq 1$, there exists
a solution to the system of polynomial equations \eqref{eq:symmetric_system_equalities}
which does not satisfy inequality \eqref{eq:symmetric_system_inequalities}.
\end{proposition}
Equipped with the definition of $\ordersym$, we now state the main result of this subsection.
\begin{theorem}  \label{theorem:convergence_sym}
Let $\pi=1/2$.
\begin{itemize}
\item[(a)] (Maximum Likelihood Estimation) We have
\begin{align*}
 \sup_{\bfeta \in H} 
\bbE_{\bfeta}\Big[\symloss_{\ordersym}(\hbfeta_n,\bfeta)\Big] \lesssim \left(\frac{\log n}{n}\right)^{\frac 1 {2\ordersym}}
 \end{align*}
where the expectation is
taken with respect to the product distribution of an i.i.d. sample $Y_{1},\ldots,Y_{n}$ 
from model \eqref{eqn:general_model}.
\item[(b)] (Minimax Lower Bound) Let $v_0 \in \text{int}(\Omega)$, $\bfeta_0 = (0,v_0, v_0)$,
and let 
$$H(\kappa) = \{\bfeta \in H: \symloss_{\ordersym}^{\ordersym}(\bfeta, \bfeta_0) \leq \kappa\}, \quad \kappa > 0.$$
Then, there exists a universal constant $c_2 > 0$ such that
\begin{align*}
\inf_{\tilde \bfeta_n}
\sup_{\substack{\bfeta \in \etaspace(c_2n^{-1/2})}} 
\bbE_{\bfeta}\Big[\symloss_{\ordersym}(\tilde\bfeta_n,\bfeta)\Big] \gtrsim \left(\frac{1}{n}\right)^{\frac 1 {2\ordersym}},
\end{align*}
where the infimum is over all sequences of estimators $\tilde \bfeta_n$ based on an i.i.d. sample $Y_1, \dots, Y_n$ from 
model \eqref{eqn:general_model}.
\end{itemize}
\end{theorem}
Similarly as in the asymmetric regime, our proof technique for Theorem 
\ref{theorem:convergence_sym}(a) hinges upon a characterization of the symmetric
loss function $\symloss_{\ordersym}$ in terms of the Total Variation
distance over the space of Gaussian mixture densities. In particular, in Theorem \ref{theorem:prelim_sym} in Appendix \ref{sec:appendix_symmetric}, we prove that
\begin{equation}
\label{eq:psi_tv_sym}
\symloss_{\ordersym}(\bfeta^{(1)}, \bfeta^{(2)}) \lesssim V\Big(g(\cdot,\bfeta^{(1)}), g(\cdot, \bfeta^{(2)})\Big)^{\frac 1 {\ordersym}}.
\end{equation}
for any $\bfeta^{(1)}, \bfeta^{(2)} \in H$. Furthermore, the bound~\eqref{eq:psi_tv_sym} is tight, which directly leads to the minimax lower bound in Theorem~\ref{theorem:convergence_sym}. Now, we obtain a specific value for the symmetric order in the following result.
\begin{proposition}
\label{proposition:symmetric_system}
Under the system of polynomial equations and inequalities \eqref{eq:symmetric_system_equalities} and
\eqref{eq:symmetric_system_inequalities}, we have $\ordersym = 4$. 
\end{proposition}
Proposition \ref{proposition:symmetric_system} together with Theorem 
\ref{theorem:convergence_sym} implies the $n^{-1/8}$ convergence rate
for the MLE under the symmetric loss function $\symloss_{\ordersym}$, up to polylogarithmic factors. 
This rate is in stark contrast to the $n^{-1/12}$ rate obtained in the asymmetric regime for a wide range of
values of $\pi$. 
In terms of parameter estimation, 
Theorem 
\ref{theorem:convergence_sym} leads to
\begin{align*}
\Big|\big|\htheta_n\big| - |\theta_{n}|\Big| \lesssim \left(\frac {\log n} n\right)^{\frac 1 8}, \quad  
 \min_{\tau \in S_2} \biggr\{|\hat v_{1,n} - v_{\tau(1),n}| + 
               |\hat v_{2,n} - v_{\tau(2),n}|
      \biggr\} \lesssim \left(\frac{\log n}{n}\right)^{\frac 1 4},
\end{align*}
with probability tending to one, where $S_2$ is the set of permutations on $\{1,2\}$.
Note that the above result merely implies a rate of convergence
for $|\widehat{\theta}_{n}|$ in absolute value. 
In the absence of absolute values, it can be shown that the above rate
becomes non-polynomial, due to the non-identifiability of the sign of $\theta_n$ under the symmetric regime.
Indeed, this situation corresponds to the use of the loss function $\assymloss_{\ordersym}$ in place
of $\symloss_{\ordersym}$
in Theorem \ref{theorem:convergence_sym}---a careful investigation of the proof reveals that 
inequality \eqref{eq:symmetric_system_inequalities} of the symmetric system would
not be needed under this loss function, which would lead to an infinite value of the symmetric order by Proposition \ref{proposition:inequality_necessity}.


In the asymmetric regime of Section \ref{sec:asymmetric}, we noted that the optimal rate of convergence
$n^{-1/12}$ under the asymmetric loss function is markedly slower than the rate $n^{-1/6}$
which is minimax optimal when the variances in model \eqref{eqn:general_model} are unknown but equal. 
Remarkably, the same behaviour does not occur in the symmetric setting:
the $n^{-1/8}$ minimax rate implied by Theorem \ref{theorem:convergence_sym}, up to polylogarithmic
factors, matches the minimax rate obtained by \cite{feller2019} when the two variances are assumed equal but unknown.
Finally, we note that the rate $(\log n/n)^{1/8}$ also matches the pointwise
rate obtained by \cite{ho2016b} when the scale parameters are not presumed equal. 
 
\section{Simulation Study}
\label{sec:simulations}
We now illustrate our theoretical results from Section \ref{sec:main} via a careful simulation study. 
We approximate the MLE using the EM algorithm
\citep{dempster1977}, tailored to the structure of model \eqref{eqn:general_model}. All simulations
below are run in Python 3.6 on a standard Linux machine. Further implementation details are relegated
to Appendix \ref{sec:appendix_numerical}.

For 100 values of the sample size $1,000 \leq n \leq 100,000$, we generate $n$
i.i.d. observations from the two-component Gaussian mixture model \eqref{eqn:general_model} with parameters
$$\theta_n = \epsilon_n (x_2^* - x_1^*), \quad 
   v_{1,n} = 1 + \epsilon_n^2 y_1^*, \quad
   v_{2,n} = 1 + \epsilon_n^2 (y_2^*+y_3^*),$$
where $s^*=(x_1^*, x_2^*, y_1^*, y_2^*, y_3^*) \in \bbR^5$ and $\epsilon_n \downarrow 0$. This setting is inspired
by our minimax analyses in Theorems \ref{theorem:convergence_asym} and \ref{theorem:convergence_sym}.
We consider two distinct settings
for $\epsilon_n$ and $s^*$. 

\begin{figure}[t]
\centering
\begin{subfigure}{0.4\textwidth}
  \centering
  \includegraphics[width=.92\linewidth]{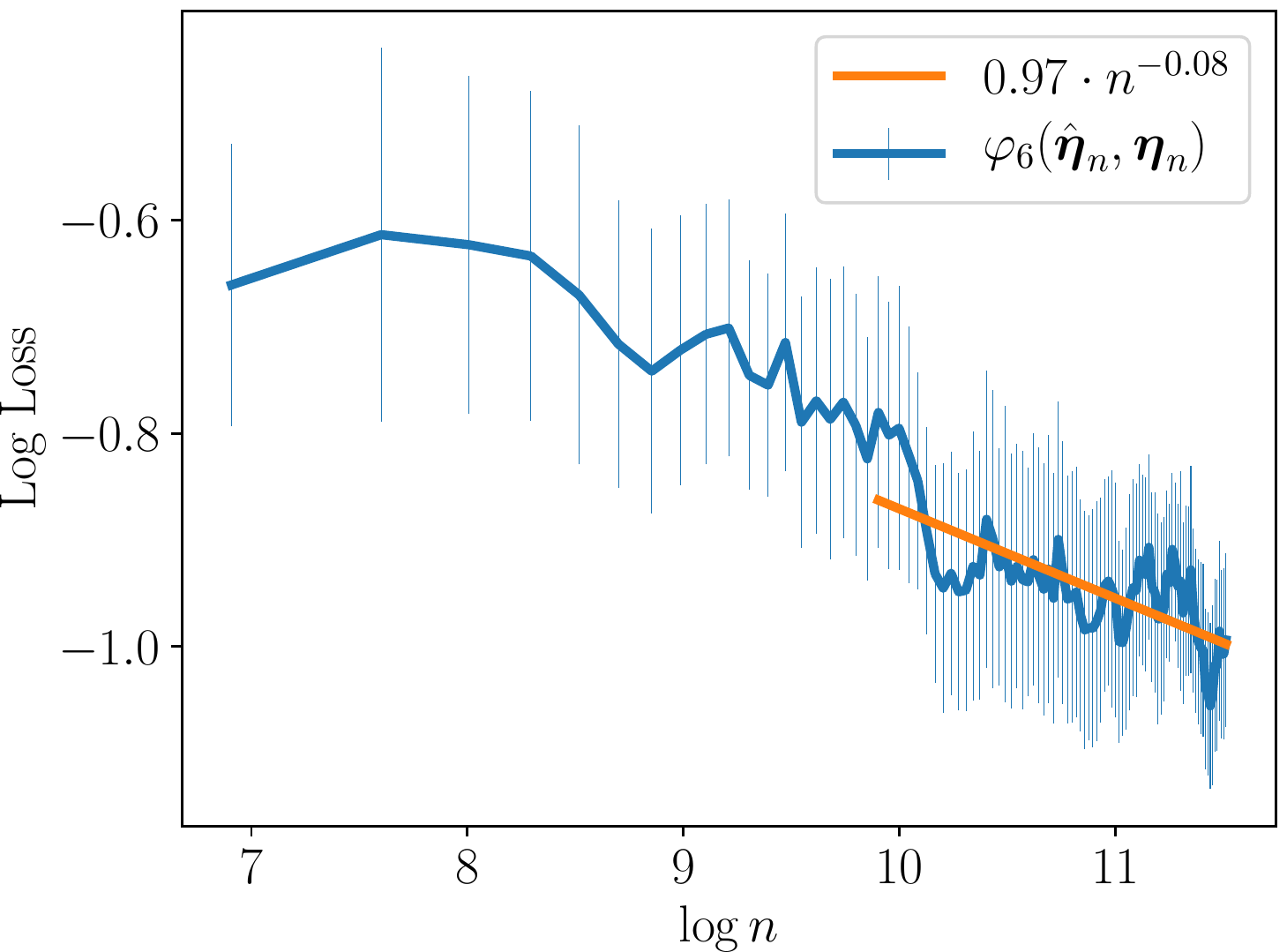}
  \caption{Model A, $\pi=0.1$.}
\end{subfigure}%
\begin{subfigure}{0.4\textwidth}
  \centering
  \includegraphics[width=.92\linewidth]{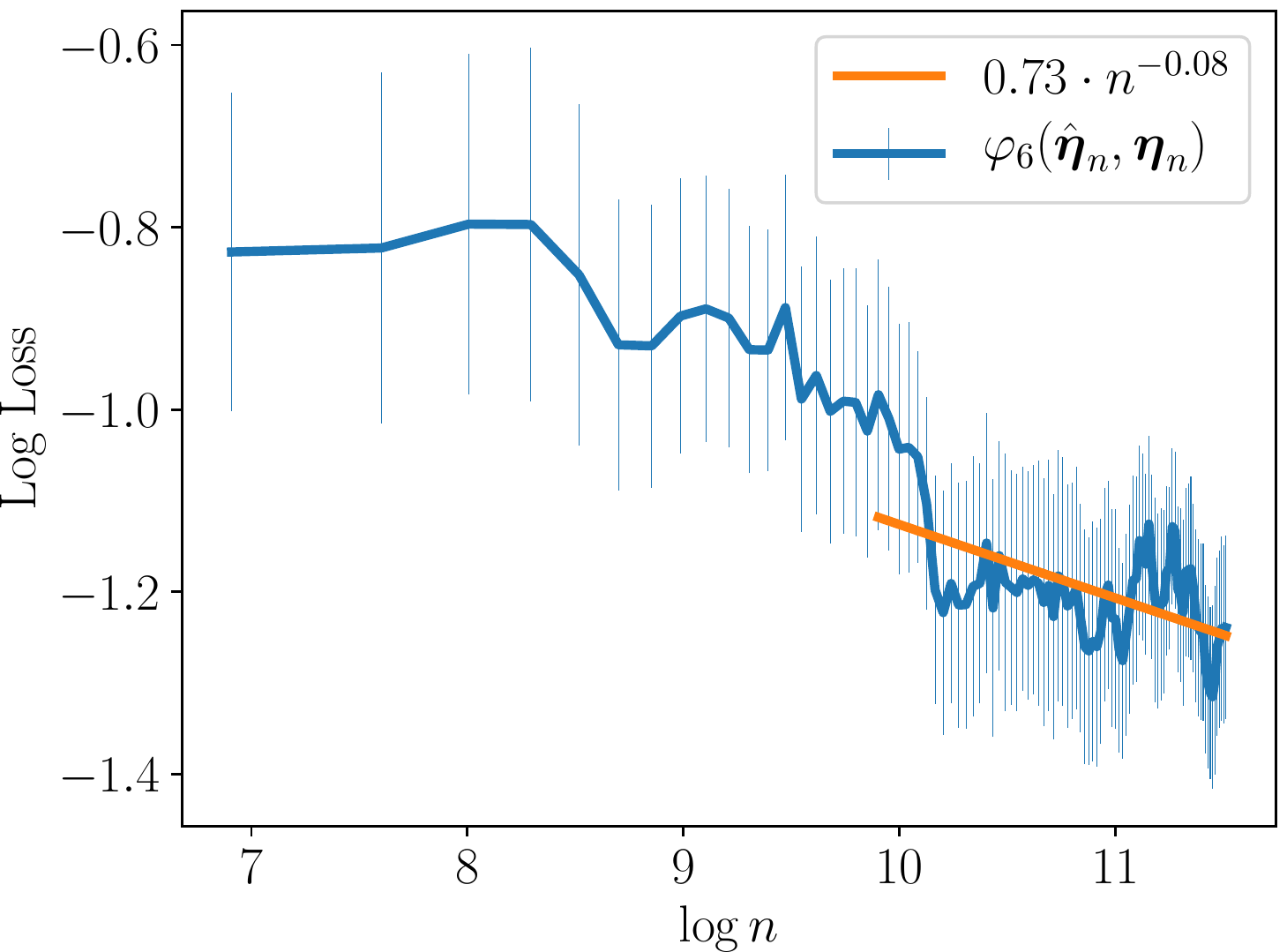}
  \caption{Model A, $\pi=0.25$.}
\end{subfigure}\par\medskip
\begin{subfigure}{0.4\textwidth}
  \centering
  \includegraphics[width=.92\linewidth]{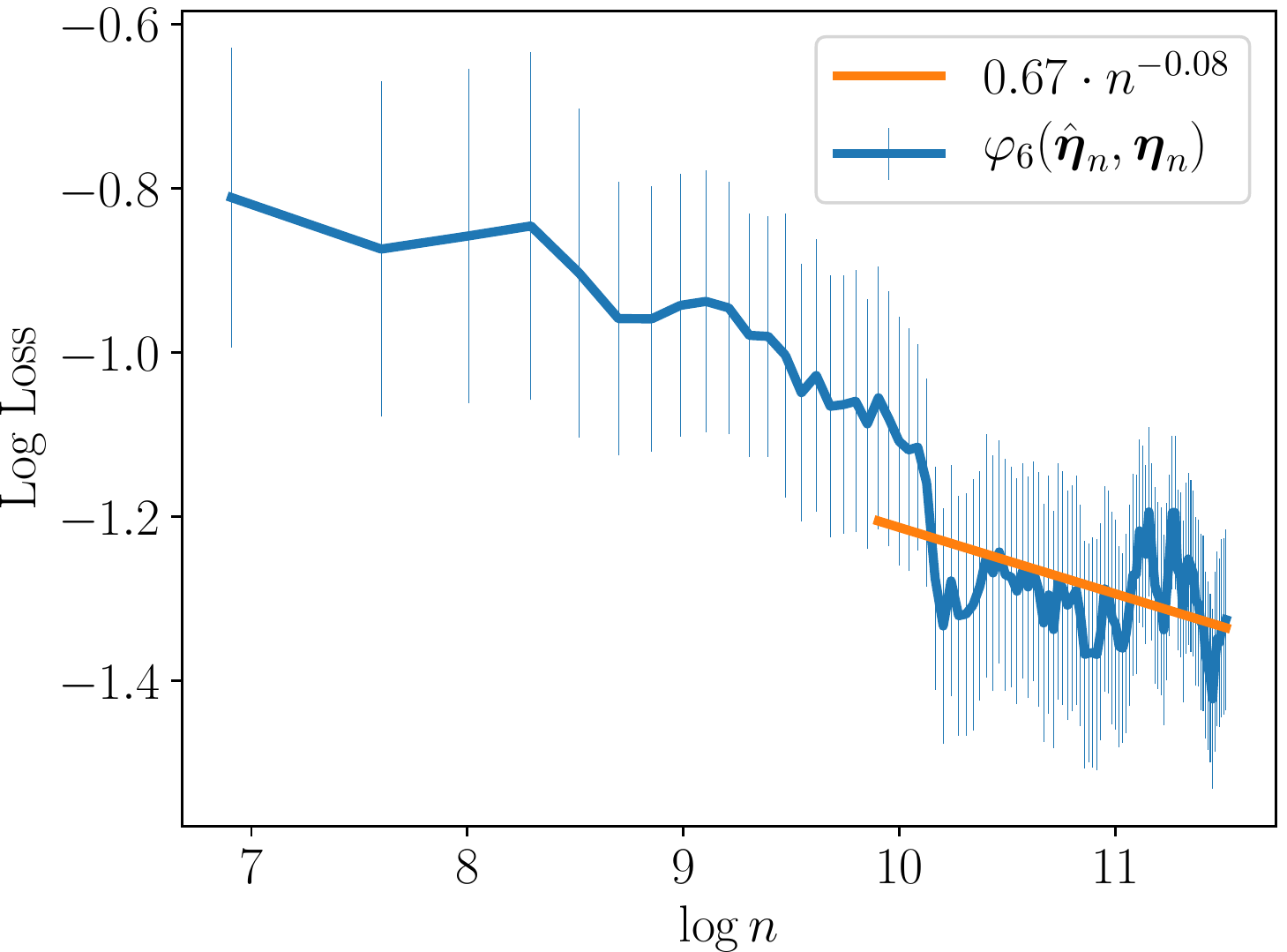}
  \caption{Model A, $\pi=0.4$.}
\end{subfigure}%
\begin{subfigure}{0.40\textwidth}
  \centering
  \includegraphics[width=.92\linewidth]{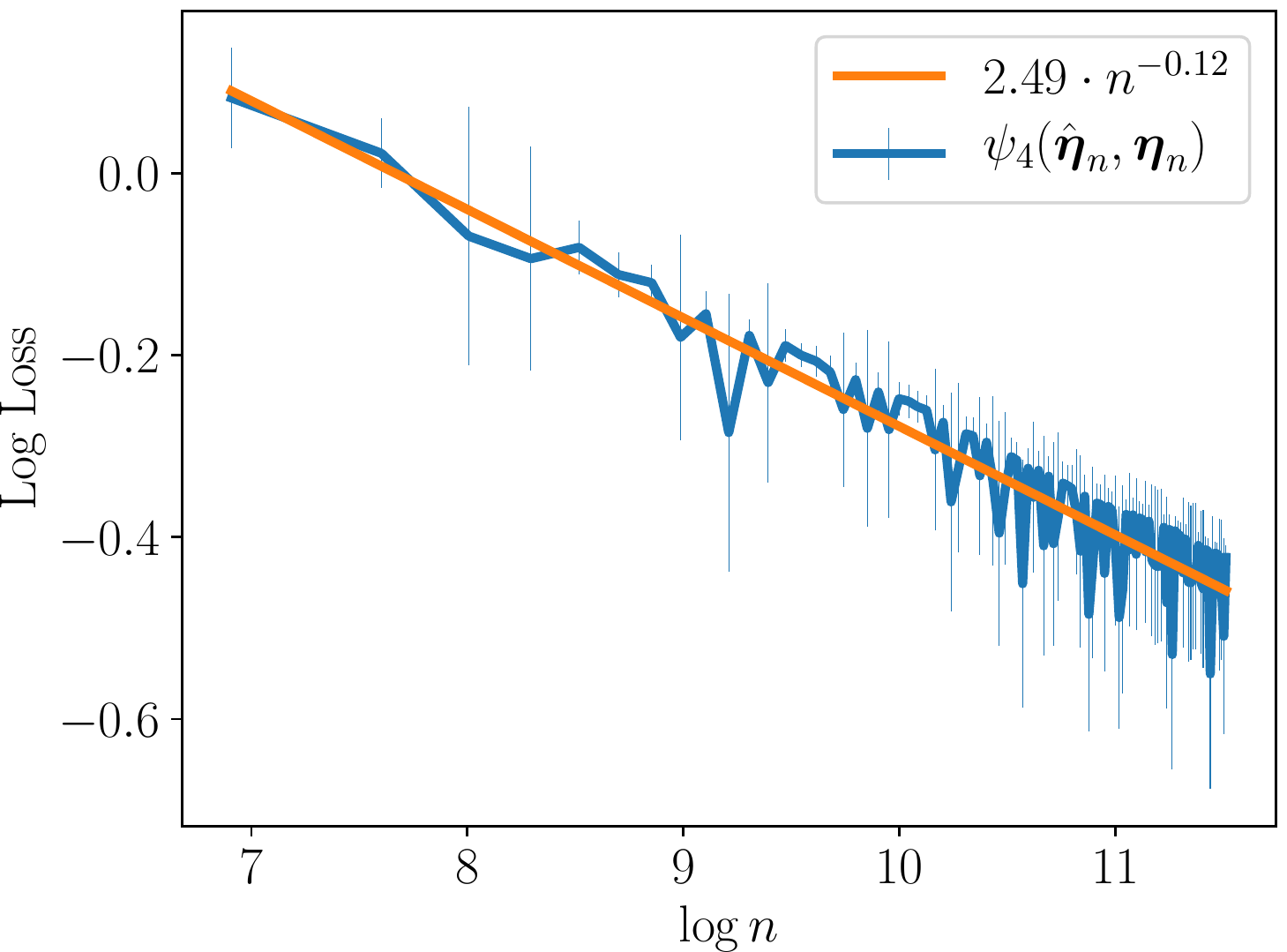}
  \caption{Model S, $\pi=0.5$.}
\end{subfigure}%
\caption{\label{fig:sims} Log-log scale plots for the simulation results under Model A ($\pi \in \{0.1, 0.25, 0.4\}$) and Model S.
For each model and sample size, the MLE $\hbfeta_n$ is computed on 10 independent samples, and its average distance from
$\bfeta_n$ is plotted with error bars representing one empirical standard deviation. 
}
\end{figure}

\textbf{Model A: Asymmetric Setting.} 
Here, we take $x_1^* = x_2^* = 0$ and $y_1^*  = -y_2^*/c, y_3^* = -y_2^*(1 + (1/c))/2, y_2^* = 0.1$.
As shown in the proof of Proposition \ref{proposition:asymmetric_system}, this choice
forms a solution to the asymmetric system of polynomial equations \eqref{eq:asymmetric_system} with $r=5$. 
Furthermore, we take $\epsilon_n = n^{-1/12}$. 
We consider three distinct
values of the mixing proportion $\pi \in \{0.1, 0.25, 0.4\}$ 
For each choice of $\pi$, we report in Figure \ref{fig:sims} (a)-(c)
the value of $\assymloss_6(\hbfeta_n,\bfeta_n)$ for all $n$ under consideration.
It can be seen that the empirical rate of convergence of the MLE is
approximately $n^{-1/12}$ for large enough $n$, under the three values of $\pi$ considered.
This rate was predicted by Theorem \ref{theorem:convergence_asym}.
While our focus in Model A is the asymmetric regime, we report in Appendix \ref{sec:appendix_numerical}
the result of this simulation with $\pi=1/2$, and we indeed observe  a markedly faster rate of convergence. 
 
\textbf{Model S: Symmetric Setting.} We now set $\pi=1/2$
and consider the setting $x_1^* = 1, x_2^* = 1.5, y_1^* = 3.5, y_2^* = 0.5, y_3=-1.5$, which
solves the symmetric system of polynomial equalities and inequalities \eqref{eq:symmetric_system_equalities}
and \eqref{eq:symmetric_system_inequalities} with $r=3$. Furthermore, we choose $\epsilon_n = n^{-1/8}$.
The empirical convergence
rate of the MLE is reported in Figure \ref{fig:sims} (d). We observe the approximate rate $n^{-1/8}$,
as anticipated by Theorem \ref{theorem:convergence_sym}.

\section{Discussion}
\label{sec:discussion}
The focus of this paper has been to derive uniform convergence rates of the maximum likelihood
estimator for parameter estimation in two-component location-scale Gaussian mixture models, 
as well as corresponding minimax lower bounds. Our analysis reveals a phase
transition in the rate of convergence depending on whether or not the mixture is symmetric. 
Specifically, we prove that the optimal rate for parameter estimation varies from $n^{-1/12}$
in the asymmetric case for a wide range of mixing proportions $\pi$, to $n^{-1/8}$ in the symmetric case, up to polylogarithmic factors. 
Key to establishing these rates is the study of certain systems of polynomial equations and inequalities, 
arising from the  dependence between the parameters of location-scale Gaussian mixtures implied
by the PDE \eqref{eq:key_pde}.

To the best of our knowledge, there are no existing works establishing 
minimax rates for parameter estimation in location-scale Gaussian mixture models with more than two components, 
except under the regime where the variances are presumed equal but unknown \citep{wu2019}.
In future work, we intend to extend the analyses of this paper to Gaussian mixtures
admitting more than two components. Furthermore, we wish to stress that the rates 
obtained in this paper are minimax, and hence are only informative about the worst-case behaviour
of parameter estimation. Mixture models are, however, notorious for admitting risk functions which
can fluctuate dramatically across the parameter space \citep{ho2019a}. 
For example, we report an extension of our simulation study in Appendix \ref{sec:appendix_numerical}
in which faster empirical rates of convergence can be observed for similar models as those of Section
\ref{sec:simulations}.  We conjecture that a more nuanced characterization
of the polynomial systems in this work would allow for 
instance-specific rates of convergence, and we are currently
exploring such directions.


\bibliography{manuscript_gaussian}

\clearpage
\appendix
\renewcommand\thefigure{\thesection.\arabic{figure}}    
{\LARGE \textbf{Appendix}\\[0.1in]}

\setcounter{figure}{0}    

In this appendix, we provide detailed proofs for all the  results in Section \ref{sec:main}. 
In Appendix \ref{sec:appendix_preliminary}, we state several existing results which will be frequently used in the sequel.
Proofs of results under the asymmetric regime can be found in Appendix \ref{sec:appendix_asymmetric}, and those
under the symmetric regime can be found in Appendix \ref{sec:appendix_symmetric}.
In Appendix \ref{sec:appendix_rasym}, we provide upper bounds on the asymmetric order $\orderassym(\pi)$
for certain values of $\pi \in(0,1/2)$. Finally, we report simulation specifications and additional
simulation results in Appendix \ref{sec:appendix_numerical}.



\section{Preliminary Results}
\label{sec:appendix_preliminary}
We begin by stating several results which will be frequently used in the sequel.
\begin{lemma}
\label{lemma:gaussian_si}
For any integer $s \geq 1$ and any $v_0 \in \Omega$, the family $\left\{\frac{\partial^\ell f}{\partial \theta^\ell}(\cdot, 0, v_0): 1 \leq \ell \leq s\right\}$
is linearly independent, in the sense that for any $\alpha_1, \dots, \alpha_s \in \bbR$, 
$$\esssup_{x \in \bbR}\left|\sum_{\ell=1}^s \alpha_\ell \frac{\partial^\ell f}{\partial \theta^\ell}(x, 0, v_0)\right| = 0
 ~~\Longrightarrow~~ \alpha_1 = \dots = \alpha_s = 0.$$
\end{lemma}
Lemma \ref{lemma:gaussian_si} follows from the strong identifiability of the location Gaussian parametric family, 
as established by \cite{chen1995} and \cite{heinrich2018}. Furthermore, we state 
a general density estimation result for the MLE $\hbfeta_n$.
\begin{lemma}
\label{lemma:hellinger_rate}There exist universal constants $c,c_1 > 0$ depending only on $\Theta,\Omega$
such that for all $u \geq c_1 \sqrt{\frac{\log n}n}$,
$$\sup_{\substack{\bfeta \in \etaspace}}\bbP\left\{h\big(g(\cdot,\hbfeta_n),
                                                                            g(\cdot,\bfeta)\big) \geq u\right\}
   \leq c\exp\left\{-\frac{nu^2}{c^2}\right\}.$$
\end{lemma}
Lemma \ref{lemma:hellinger_rate} may be obtained by combining a guarantee for the maximum likelihood
density estimator (see for instance Theorem 7.4 of \cite{vandegeer2000}) with bracketing numbers for 
classes of Gaussian mixture densities \citep{genovese2000, ghosal2001}.
See Theorem 4.1 of \cite{ho2016c} for further details.

\section{Proofs under the Asymmetric Regime}
\label{sec:appendix_asymmetric}
\subsection{Proof of Theorem \ref{theorem:convergence_asym} }
The essence of Theorem \ref{theorem:convergence_asym}.(a) is contained in the following result. 
\begin{theorem}
\label{theorem:prelim_asym}
Let $\pi\in (0,1/2)$. Then,
\begin{eqnarray}
\inf \limits_{\bfeta^{(1)}, \bfeta^{(2)} \in \etaspace} \dfrac{V\parenth{g(\cdot,\bfeta^{(1)}),g(\cdot,\bfeta^{(2)})}}{
\assymloss_{\orderassym(\pi)}^{\orderassym(\pi)}(\bfeta^{(1)}, \bfeta^{(2)})} > 0. \nonumber
\end{eqnarray}
\end{theorem}

We begin by proving Theorem \ref{theorem:prelim_asym}, and we then prove Theorem \ref{theorem:convergence_asym} below.

\paragraph{PROOF OF THEOREM \ref{theorem:prelim_asym}} 

We write $\rbar = \orderassym(\pi)$ throughout the proof. Assume by way of a contradiction 
that the claim does not hold. 
We may then find sequences $\bfeta_n^{(1)} = (\theta_n^{(1)}, v_{1,n}^{(1)}, v_{2,n}^{(1)}) ,
\bfeta_n^{(2)} = (\theta_n^{(2)}, v_{1,n}^{(2)}, v_{2,n}^{(2)}) \in \etaspace$, $n\geq 1$, such that
\begin{eqnarray}
\label{eq:non_symmetric_a1_assm}
\frac 1 {D_n} V\parenth{g(\cdot,\bfeta_n^{(1)}), g(\cdot,\bfeta_n^{(2)})} \to 0 
\end{eqnarray}
as $n \to \infty$, where $D_n = \assymloss_{\rbar}^{\rbar}(\bfeta_n^{(1)}, \bfeta_n^{(2)})$.
For the convenience of presentation, we only consider the most challenging setting $\theta_{n}^{(j)} \to 0$, $v_{i,n}^{(j)} \to v_{0}$
for all $i,j = 1,2$,  for some $v_{0} \in \Omega$.

Our proof will rely on the following setup. By Taylor expansion up to order $\rbar$, we have
for all $x \in \bbR$,
\begin{align}
\nonumber
g(x,  &\bfeta_n^{(1)}) - g(x, \bfeta_n^{(2)}) \\
\nonumber
 &=    \pi \Big\{ f(x, -\theta_n^{(1)}, v_{1,n}^{(1)}) - f(x, -\theta_n^{(2)}, v_{1,n}^{(2)})\Big\} + 
    (1-\pi)\Big\{ f(x, c\theta_n^{(1)}, v_{2,n}^{(1)}) - f(x, c\theta_n^{(2)},  v_{2, n}^{(2)})\Big\} \\
\nonumber
 &= \pi \sum_{|\alpha| =1}^{\rbar} \frac{(\theta_n^{(2)} - \theta_n^{(1)})^{\alpha_1} (v_{1,n}^{(1)} - v_{1,n}^{(2)})^{\alpha_2}}{\alpha_1!\alpha_2!}
      \frac{\partial^{|\alpha| f}}{\partial \theta^{\alpha_1} \partial v^{\alpha_2}}(x, -\theta_n^{(2)}, v_{n,1}^{(2)}) + R_{1,n}(x)\\
\nonumber
 &+ (1-\pi) \sum_{|\alpha| =1}^{\rbar} \frac{c^{\alpha_1}(\theta_n^{(1)} - \theta_n^{(2)})^{\alpha_1}(v_{2,n}^{(1)} - v_{2,n}^{(2)})^{\alpha_2}}{
 			\alpha_1!\alpha_2!} \frac{\partial^{|\alpha|} f}{\partial \theta^{\alpha_1}\partial v^{\alpha_2}} (x, c\theta_n^{(2)}, v_n^{(2)}) + R_{2,n}(x)\\
\nonumber 			
 &= \pi \sum_{|\alpha| =1}^{\rbar} \frac 1 {2^{\alpha_2}\alpha_1!\alpha_2!} (\theta_n^{(2)} - \theta_n^{(1)})^{\alpha_1}(v_{1,n}^{(1)} - v_{1,n}^{(2)})^{\alpha_2}
  		\frac{\partial^{\alpha_1 + 2\alpha_2} f }{\partial\theta^{\alpha_1 + 2\alpha_2}} (x, -\theta_n^{(2)}, v_{1,n}^{(2)}) + R_{1,n}(x) \\ 
  		&+ 
 	 (1-\pi) \sum_{|\alpha| =1}^{\rbar} 
 	   \frac 1 {2^{\alpha_2}\alpha_1!\alpha_2!} (\theta_n^{(1)} - \theta_n^{(2)})^{\alpha_1} (v_{2,n}^{(1)} - v_{2,n}^{(1)})^{2\alpha_2}
 	   \frac{\partial f^{\alpha_1 + 2\alpha_2}}{\partial\theta^{\alpha_1+2\alpha_2}}(x, c\theta_n^{(2)}, v_{2,n}^{(2)})+R_{2,n}(x),
\label{eq:pf_asym_expansion} 	   
\end{align}
where the last equality is due to the key PDE \eqref{eq:key_pde} for location-scale Gaussian densities.
Furthermore, $R_{1,n}, R_{2,n}$ are Taylor remainders satisfying
$$\max \big\{\norm {R_{1,n}}_\infty, \norm{R_{2,n}}_\infty\big\} = O\Big(|\theta_n^{(1)} - \theta_n^{(2)}|^{\rbar+\gamma} +
                                                                 |v_{1,n}^{(1)} - v_{1,n}^{(2)}|^{\rbar+\gamma} + 
                                                                 |v_{2,n}^{(2)} - v_{2,n}^{(2)}|^{\rbar+\gamma}\Big),$$
for some $\gamma > 0$. 
Now, by a further Taylor expansion to order $\rbar-|\alpha|$ of the partial derivatives
appearing in \eqref{eq:pf_asym_expansion}, we also have 
for all $x \in \bbR$ and all $1 \leq |\alpha| \leq \rbar,$
\begin{align}
\label{eq:pf_asym_partial_expansion}
&\dfrac{\partial^{\alpha_{1}+2\alpha_{2}}{f}}{\partial{\theta^{\alpha_{1}+2\alpha_{2}}}}(x,c\theta_{n}^{(2)},v_{2,n}^{(2)}) \nonumber \\
&\quad= \sum \limits_{|\beta|=0}^{\rbar-|\alpha|}
		\dfrac{(c+1)^{\beta_{1}}(\theta_{n}^{(2)})^{\beta_{1}}(v_{2,n}^{(2)} - v_{1,n}^{(2)})^{\beta_{2}}}{2^{\beta_{2}}\beta_{1}!\beta_{2}!}\dfrac{\partial^{\alpha_{1}+2\alpha_{2}+\beta_{1}+2\beta_{2}}{f}}{\partial{\theta^{\alpha_{1}+2\alpha_{2}}}}(x,-\theta_{n}^{(2)},v_{1,n}^{(2)}) 
		+ R_{2, n,\alpha}(x),
\end{align}
where we have again used the PDE \eqref{eq:key_pde}. Here, 
$R_{2,n,\alpha}$ are Taylor remainders such that 
$\|R_{2,n,\alpha}\|_{\infty} = O\parenth{|\theta_{n}^{(2)}|^{r-|\alpha|+\gamma}+|v_{2,n}^{(2)} - v_{1,n}^{(2)}|^{r-|\alpha|+\gamma}}$ for all $1 \leq |\alpha| \leq r$. 
Combining the expansions in \eqref{eq:pf_asym_expansion} and \eqref{eq:pf_asym_partial_expansion}, 
we obtain the representation
$$\hspace{ - 2 em} \dfrac{g(x,\bfeta_n^{(1)}) - g(x,\bfeta_n^{(2)})}{D_{n}} 
 = \frac 1 {D_n} \left[\sum \limits_{\ell=1}^{2\rbar} A_{n,\ell}\dfrac{\partial^{\ell}{f}}{\partial{\theta^{\ell}}}(x,-\theta_{n}^{(2)},v_{1,n}^{(2)})+R_n(x)\right],
$$
where the formulations of $A_{n,\ell}$ and $R_n(x)$ are as follows
\begin{eqnarray}
& & \hspace{-2 em} A_{n,\ell}= \pi \sum \limits_{\alpha_{1},\alpha_{2}} \dfrac{(\theta_{n}^{(2)} - \theta_{n}^{(1)})^{\alpha_{1}}(v_{1,n}^{(1)} - v_{1,n}^{(2)})^{\alpha_{2}}}{2^{\alpha_{2}}\alpha_{1}!\alpha_{2}!} \nonumber \\
& & +(1- \pi) \sum \limits_{\alpha_{1},\alpha_{2},\beta_{1},\beta_{2}}\dfrac{c^{\alpha_{1}}(c+1)^{\beta_{1}}(\theta_{n}^{(1)} - \theta_{n}^{(2)})^{\alpha_{1}}(\theta_{n}^{(2)})^{\beta_{1}}(v_{2,n}^{(1)} - v_{2,n}^{(2)})^{\alpha_{2}}(v_{2,n}^{(2)}  - v_{1,n}^{(2)})^{\beta_{2}}}{2^{\alpha_{2}+\beta_{2}}\alpha_{1}!\alpha_{2}!\beta_{1}!\beta_{2}!} \nonumber \\
& & \hspace{-2 em} R_n(x) = \pi R_{1,n}(x) + (1-\pi)R_{2,n}(x) + \sum \limits_{|\alpha| \leq \rbar}\dfrac{1}{2^{\alpha_{2}}} \dfrac{c^{\alpha_{1}}(\theta_{n}^{(1)} - \theta_{n}^{(2)})^{\alpha_{1}}(v_{2,n}^{(1)} - v_{2,n}^{(2)})^{\alpha_{2}}}{\alpha_{1}!\alpha_{2}!}R_{2,n,\alpha}(x) \nonumber
\end{eqnarray}
for all $1 \leq \ell \leq 2\rbar$ where the ranges of $\alpha_{1},\alpha_{2}$ in the first sum of $A_{n,\ell}$ satisfy 
$\alpha_{1}+2\alpha_{2}=\ell$, 
$1 \leq \alpha_{1}+\alpha_{2} \leq \rbar$ and 
the ranges of $\alpha_{1},\alpha_{2},\beta_{1},\beta_{2}$ in the second sum of $A_{n,\ell}$ 
satisfy $\alpha_{1}+\beta_{1}+2\alpha_{2}+2\beta_{2} = \ell$, $1 \leq \alpha_{1}+\alpha_{2} \leq \rbar$, 
and $0 \leq \beta_{1}+\beta_{2} \leq \rbar - (\alpha_{1}+\alpha_{2})$. 
 
We now prove the claim by considering 
the following settings regarding
$\theta_{n}^{(1)}, \theta_{n}^{(2)}, v_{1,n}^{(1)}, v_{1,n}^{(2)}, v_{2,n}^{(1)}$, and $v_{2,n}^{(2)}$.

\paragraph{Case a.1:} $\theta_{n}^{(2)}/\theta_{n}^{(1)} \not \to 1$ and $|v_{2,n}^{(2)} - v_{1,n}^{(2)}|/\max \left\{|v_{1,n}^{(1)} - v_{1,n}^{(2)}|, |v_{2,n}^{(1)} - v_{2,n}^{(2)}|\right\} \not \to \infty$ as $n \to \infty$. 

In this case, it is a straightforward verification
 that $\|R_n\|_{\infty}/D_{n} \to 0$. 
We further claim that there exists $1 \leq \ell \leq \rbar$ such that
$A_{n,\ell}/D_{n} \not\to 0$. 
Assume by way of a contradiction that $A_{n,\ell}/D_{n}\to 0$ for all such $l$. We denote
\begin{eqnarray}
\overline{M}_{n} = \max\left\{|\theta_{n}^{(2)} - \theta_{n}^{(1)}|, |v_{1,n}^{(1)} - v_{1,n}^{(2)}|^{1/2}, |v_{2,n}^{(1)} - v_{2,n}^{(2)}|^{1/2}\right\}. \nonumber
\end{eqnarray}
From the assumption of Case a.1, 
we have $|v_{2,n}^{(2)} - v_{1,n}^{(2)}|/\overline{M}_{n}^{2} \not \to \infty$ and $|\theta_{n}^{(2)}|/\overline{M}_{n} \not \to \infty$. 
Therefore, there exist $x_1, x_2, y_1, y_2, y_3 \in \bbR$ such that
$$(\theta_{n}^{(2)} - \theta_{n}^{(1)})/\overline{M}_{n} \to x_{1}, \quad 
   \theta_{n}^{(2)}/\overline{M}_{n} \to x_{2}, $$
and,   
$$ (v_{1,n}^{(1)} - v_{1,n}^{(2)})/\overline{M}_{n}^{2} \to y_{1},\quad 
   (v_{2,n}^{(1)} - v_{2,n}^{(2)})/\overline{M}_{n}^{2} \to y_{2}, \quad
   (v_{2,n}^{(2)} - v_{1,n}^{(2)})/\overline{M}_{n}^{2} \to y_{3}.$$ 
From the definition of $\overline{M}_{n}$, at least one of $x_{1},y_{1},y_{2}$ is different from 0. Now, by dividing both the numerator and the denominator of $A_{n,\ell}$ ($1 \leq \ell \leq \rbar$) by $\overline{M}_{n}^{\ell}$, and using the fact that
$A_{n,\ell}/D_{n}\to 0$ for all $1 \leq \ell \leq \rbar$, we obtain the asymmetric system of polynomial equations
\begin{eqnarray}
\hspace{ - 3 em} (1-\pi) \sum \limits_{\alpha_{1},\alpha_{2},\beta_{1},\beta_{2}} \dfrac{1}{2^{\alpha_{2}+\beta_{2}}} \dfrac{c^{\alpha_{1}}(c+1)^{\beta_{1}}(-x_{1})^{\alpha_{1}}x_{2}^{\beta_{1}}y_{2}^{\alpha_{2}}y_{3}^{\beta_{2}}}{\alpha_{1}!\alpha_{2}!\beta_{1}!\beta_{2}!} \nonumber \\
+ \pi \sum \limits_{\alpha_{1},\alpha_{2}} \dfrac{1}{2^{\alpha_{2}}}\dfrac{x_{1}^{\alpha_{1}}y_{1}^{\alpha_{2}}}{\alpha_{1}!\alpha_{2}!}  = 0 \label{eqn:system_polynomial_non_symmetric_case}
\end{eqnarray}
for all $1 \leq \ell \leq \rbar$. 
By the choice of $\rbar$, 
this system does not admit any non-trivial solutions, thus
$x_{1}=y_{1}=y_{2}=0$, contradicting the fact that at least one among $x_{1}, y_{1}, y_{2}$ is different from 0. 
Therefore, not all the coefficients $A_{n,\ell}/D_{n}$  tend to 0 as $n \to \infty$. 

Letting $m_{n} = D_{n}/\max \limits_{1 \leq \ell \leq 2\rbar} |A_{n,\ell}|$, it follows that $m_{n} \not \to \infty$. Now, for all $x \in \bbR$, we have that
\begin{eqnarray}
\label{eq:non_symmetric_a1_step1}
\frac{m_{n}}{D_n} \sum_{\ell=1}^{2\rbar} A_{n,\ell}\dfrac{\partial^{\ell}{f}}{\partial{\theta^{\ell}}}(x,-\theta_{n}^{(2)},v_{1,n}^{(2)}) 
\to \sum \limits_{\ell=1}^{2\rbar} \tau_{\ell}\dfrac{\partial^{\ell}{f}}{\partial{\theta^{\ell}}}(x,0,v_{0}) 
\end{eqnarray}
for some coefficients $\tau_{\ell} \in \bbR$ which are not all 0. On the other hand, since $m_n \not\to\infty$,
assumption \eqref{eq:non_symmetric_a1_assm} implies
\begin{align}
\nonumber
\frac{m_n}{D_n} V\Big(g(\cdot,&\bfeta_n^{(1)}), g(\cdot,\bfeta_n^{(2)})\Big) \\ 
\nonumber
 &= \int \frac{m_n}{D_n} \left|g(x,\bfeta_n^{(1)}) - g(x,\bfeta_n^{(2)})\right|dx  \\
 &= \int \frac{m_n}{D_n} \left|\sum_{\ell=1}^{2\rbar} A_{n,\ell}\dfrac{\partial^{\ell}{f}}{\partial{\theta^{\ell}}}(x,-\theta_{n}^{(2)},v_{1,n}^{(2)})
  + R_n(x)\right|dx  \to 0.
 \label{eq:non_symmetric_a1_step2}
\end{align}
By Fatou's Lemma, the integrand in \eqref{eq:non_symmetric_a1_step2} 
vanishes to zero almost everywhere, and since $\norm{R_n}_\infty /D_n \to 0$, we arrive at
$$\frac{m_{n}}{D_n} \left|\sum_{\ell=1}^{2\rbar} A_{n,\ell}\dfrac{\partial^{\ell}{f}}{\partial{\theta^{\ell}}}(x,-\theta_{n}^{(2)},v_{1,n}^{(2)})\right|
 \to 0,$$
 for almost every $x \in \bbR$.
Combining this fact with \eqref{eq:non_symmetric_a1_step1} then yields
 \begin{align}
\sum \limits_{\ell=1}^{2\rbar} \tau_{\ell}\dfrac{\partial^{\ell}{f}}{\partial{\theta^{\ell}}}(x,0,v_{0}) = 0. \nonumber
\end{align} 
However, by Lemma \ref{lemma:gaussian_si} we have $\tau_{\ell} = 0$ for all
$1 \leq \ell \leq 2\rbar$, which is a contradiction. Therefore, Case a.1 does not hold. 
\paragraph{Case a.2:} $\theta_{n}^{(2)}/\theta_{n}^{(2)} \not \to 1$ and $|v_{2,n}^{(2)} - v_{1,n}^{(2)}|/\max \left\{|v_{1,n}^{(1)} - v_{1,n}^{(2)}|, |v_{2,n}^{(1)} - v_{2,n}^{(2)}|\right\} \to \infty$ as $n \to \infty$. 

Unlike Case a.1, we will not necessarily have $\|R_n\|_{\infty}/D_{n} \to 0$ in this case. 
Our approach instead hinges upon the following Lemma.
\begin{lemma}
\label{lemma:case_a2_asym}
Under Case a.2, we have
$$\max \limits_{1 \leq \ell \leq 2\rbar} |A_{n,\ell}|/D_{n} \not \to 0,\quad \text{and}, \quad  
\|R_n\|_{\infty}/\max \limits_{1 \leq \ell \leq 2\rbar} |A_{n,\ell}| \to 0.$$
\end{lemma}
The proof of Lemma \ref{lemma:case_a2_asym} appears in Appendix \ref{sec:lemmas_asym} below. 
Writing $m_n' = D_n / \max \limits_{1 \leq \ell \leq 2\rbar}|A_{n,\ell}|$, 
Lemma \ref{lemma:case_a2_asym} implies $m_n' \not\to \infty$ and $m_n' \norm{R_n}/D_n \to 0$, 
thus following similar steps as in Case a.1, we arrive at
$$\frac{m_n'}{D_n} \left[g(x, \bfeta_n^{(1)}) - 
                        g(x, \bfeta_n^{(2)})\right] \to 
                        \sum_{\ell=1}^{2\rbar} \tau_\ell' \frac{\partial^\ell f}{\partial \theta^\ell}(x, 0, v_0),$$
for some $\tau_\ell' \in \bbR$ not all zero. Then, similarly as in Case a.1, 
Fatou's Lemma combined with the hypothesis \eqref{eq:non_symmetric_a1_assm}
implies
$$\sum_{\ell=1}^{2\rbar} \tau_\ell' \frac{\partial f}{\partial\theta^\ell} (x, 0, v) = 0,$$
for almost every $x \in \bbR$, 
contradicting the result of Lemma \ref{lemma:gaussian_si}.
Thus, Case a.2 cannot hold.

\paragraph{Case a.3:} $\theta_{n}^{(2)}/\theta_{n}^{(1)} \to 1$ as $n \to \infty$. Similarly 
to Case a.2,  $\|R \|_{\infty}/D_{n}$ does not generically tend to 0 in this case. We prove the following Lemma in
Appendix \ref{sec:lemmas_asym} below. 
\begin{lemma}
\label{lemma:case_a3_asym}
Under Case a.3, we have
$$\max \limits_{1 \leq \ell \leq 2\rbar} |A_{n,\ell}|/D_{n} \not \to 0, \quad \text{and} \quad 
\|R_n\|_{\infty}/\max \limits_{1 \leq \ell \leq 2\rbar} |A_{n,\ell}| \to 0.$$
\end{lemma}
By the same argument as Case a.2, it can readily be shown that Case a.3 does not hold. We have thus derived a contradiction
with \eqref{eq:non_symmetric_a1_assm}. The claim follows. 
\hfill$\square$\\


We are now in a position to prove the main result of this section.

\paragraph{PROOF OF THEOREM \ref{theorem:convergence_asym}} 

Throughout the proof, we write $\rbar = \orderassym(\pi)$.

(a) Let $\bfeta \in \etaspace$.
Theorem \ref{theorem:prelim_asym} implies the existence of a universal constant $C_1 > 0$, depending only on $\Omega, \Theta$, such that
the MLE $\hbfeta_n$ of $\bfeta$ satisfies
\begin{align*}
\assymloss_{\rbar}^{\rbar}(\hbfeta_n, \bfeta)
   \leq C_1 V\big(g(\cdot,\hbfeta_n),g(\cdot,\bfeta)\big)
   \leq C_1 h\big(g(\cdot,\hbfeta_n),g(\cdot,\bfeta)\big),
\end{align*}
where the last inequality of the above display is due to the well-known inequality $V \leq  h$. 
Invoking Lemma \ref{lemma:hellinger_rate}, we obtain for all $u \geq c_1 \sqrt{\log n/n}$,
$$\assymloss_{\rbar}(\hbfeta_n, \bfeta)  \leq  (C_1 u)^{\frac 1 {\rbar}} $$
with probability at least $1-c \exp(-2nu^2/c^2)$. Integrating this tail probability
inequality to obtain a bound in expectation
readily yields the claim.

(b) By definition of $\rbar$, there exists a non-trivial solution
$(x_1^*, x_2^*, y_1^*, y_2^*, y_3^*) \in \bbR^5$  to the system
of polynomials \eqref{eq:asymmetric_system} with respect to the choice $r = \rbar-1$. 
Set $\bfeta_n^{(i)} = (\theta_n^{(i)}, v_{1,n}^{(i)}, v_{2,n}^{(i)}) \in H$
for $i=1,2$, where, for $\epsilon_n = n^{-1/2\rbar}$,
$$\theta_n^{(1)} = \epsilon_n (x_2^*-x_1^*),~~
  \theta_n^{(2)} = \epsilon_n x_2^*,$$
and where, recalling that $\bfeta_0=(0,0,v_0)$, 
$$ v_{1,n}^{(1)} = \epsilon_n^2 y_1^* + v_0 , ~~
   v_{2,n}^{(1)} = \epsilon_n^2 (y_2^* + y_3^*) + v_0 , ~~
   v_{1,n}^{(2)} = v_0 , ~~
   v_{2,n}^{(2)} = \epsilon_n^2 y_3^* + v_0.$$
The definition of $\epsilon_n$ then implies
that $\assymloss_{\rbar}(\bfeta_n^{(j)},\bfeta_0) \leq c_1 n^{-1/2}$ for some $c_1 > 0$, for $j=1,2$. 
Furthermore, this choice of parameters leads to the following key equalities
which will be used in the sequel
\begin{equation}
\label{eq:pf_minimax_asym_identity} 
\theta_n^{(2)} - \theta_n^{(1)} = \epsilon_n x_1^*, \quad
  v_{1,n}^{(2)} - v_{1,n}^{(2)} = \epsilon_n^2 y_1^*, \quad
  v_{2,n}^{(1)} - v_{2,n}^{(2)} = \epsilon_n^2 y_2^*,\quad
  v_{2,n}^{(2)} - v_{1,n}^{(2)} = \epsilon_n^2 y_3^*.
\end{equation}
Now, it follows from Le Cam's Inequality (\cite{tsybakov2008}, Theorem 2.2) 
that 
\begin{align*}
\inf_{\hat \bfeta_n} &\sup_{\substack{\bfeta \in \etaspace \\ \assymloss_{\rbar}(\bfeta,\bfeta_0)\leq c_1n^{-1/2}}}
 \bbE_{\bfeta} \big[\assymloss_{\rbar}(\hbfeta_n, \bfeta) \big]
 \geq \assymloss_{\rbar}(\bfeta_n^{(1)}, \bfeta_n^{(2)}) \Big(1 - V\big(g(\cdot, \bfeta^{(1)}_n), g(\cdot, \bfeta_n^{(2)})\big)\Big)
\end{align*}
Furthermore, notice that the identities \eqref{eq:pf_minimax_asym_identity} imply
$$\assymloss_{\rbar}^{\rbar}(\bfeta_n^{(1)}, \bfeta_n^{(2)}) = 
\epsilon_n^{\rbar} \Big[ |x_1^*|^{\rbar} + |y_2^*|^{\rbar/2} + |y_3^*|^{\rbar/2}\Big].$$
Since $(x_1^*, x_2^*, y_1^*, y_2^*, y_3^*)$ form a non-trivial solution, the factor
in brackets in the above display is nonzero, hence 
$\assymloss_{\rbar}(\bfeta_n^{(1)}, \bfeta_n^{(2)}) \asymp \epsilon_n$.
This fact combined
with the inequality $V \leq h$ and the tensorization property of the Hellinger distance (\cite{tsybakov2008}, p. 83)  imply
\begin{align*}
\inf_{\hat \bfeta_n} \sup_{\substack{\bfeta \in \etaspace \\ \assymloss_{\rbar}(\bfeta,\bfeta_0)\leq c_1n^{-1/2}}}& \bbE_{\bfeta}
 \big[\assymloss_{\rbar}(\hbfeta_n, \bfeta) \big]\\
 &\gtrsim  n^{-1/{2\rbar}}
 \left( 1 - \sqrt{1-\left[1-h^2\Big(g(\cdot,\bfeta_n^{(1)}), g(\cdot,\bfeta_n^{(2)})\Big)\right]^n}\right).
 \end{align*}
Notice that the right-hand side of the above display will be of order $n^{-1/2 \rbar}$ provided 
\begin{equation}
\label{eq:hellinger_goal}
h^2\Big(g(\cdot,\bfeta_n^{(1)}), g(\cdot,\bfeta_n^{(2)})\Big) \lesssim \frac 1 n.
\end{equation}
To prove the claim, it will therefore suffice to prove \eqref{eq:hellinger_goal}. 
Notice that
\begin{equation}
h^2\Big(g(\cdot,\bfeta_n^{(1)}), 
 g(\cdot,\bfeta_n^{(2)})\Big)
 = \int \frac{\left[g(x,\bfeta_n^{(1)}) -
 g(x, \bfeta_n^{(2)})\right]^2}{
 {\left[\sqrt{g(x,\bfeta_n^{(1)})} +
  \sqrt{g(x,\bfeta_n^{(2)})}\right]^2 }} dx
 \end{equation}
We begin by analyzing the numerator of the integrand in the above display.
By a similar Taylor expansion as in the proof of Theorem  \ref{theorem:prelim_asym},
but now up to order $\rbar -1$, we have for any $x \in \bbR$,
\begin{align*}
g(x, \bfeta_n^{(1)}) - g(x, \bfeta_n^{(2)}) 
 &= \sum_{\ell=1}^{2 \rbar-2} A_{n,\ell} \frac{\partial^\ell f}{\partial\theta^\ell} (x, -\theta_n^{(2)}, v_{1,n}^{(2)}) + R_n(x),
\end{align*}
where the Taylor remainder $R_n$ is given by
$$R_n(x) = \pi R_{1,n}(x) + (1-\pi)R_{2,n}(x) +
\sum \limits_{|\alpha| \leq \rbar}  \dfrac{c^{\alpha_{1}}(\theta_{n}^{(1)} - \theta_{n}^{(2)})^{\alpha_{1}}(v_{2,n}^{(1)} - v_{2,n}^{(2)})^{\alpha_{2}}}{2^{\alpha_{2}}\alpha_{1}!\alpha_{2}!}R_{2,n,\alpha}(x),$$
and where
\begin{align*}
R_{j,n}(x) &= \sum_{|\beta|=\rbar}\frac{(\theta_n^{(2)} - \theta_n^{(1)})^{\beta_1}(v_{j,n}^{(1)}-v_{j,n}^{(2)})^{\beta_2}}{\beta_1!\beta_2!} 
 \\ &\times
					\int_0^1 (1-t)^{\rbar-1} 
					\frac{\partial^{\rbar} f}{\partial\theta^{\beta_1} \partial v^{\beta_2}}\left(x, -\theta_n^{(2)}+ t(\theta_n^{(1)}-\theta_n^{(2)}), v_{1,n}^{(2)} + t(v_{j,n}^{(1)}-v_{j,n}^{(2)})\right) dt,\quad j=1,2,\\
R_{2,n,\alpha}(x) &=  \frac{(c+1)^{\rbar-|\alpha|}(\theta_n^{(2)})^{\rbar-|\alpha|}}{(\rbar-|\alpha|)!}  \\ &\times
					\int_0^1 (1-t)^{\rbar-1} \frac{\partial^{\rbar} f}{\partial\theta^{\rbar-\alpha_2} \partial v^{\alpha_2}}
					 \Big(x, -\theta_n^{(2)}+ t(\theta_n^{(1)}-\theta_n^{(2)}), v_{2,n}^{(2)} + t(v_{2,n}^{(1)}-v_{2,n}^{(2)})\Big) dt.		 
\end{align*}
Furthermore, the coefficients $A_{n,\ell}$ are given by
\begin{eqnarray}
& & \hspace{-2 em} A_{n,\ell}= \pi \sum \limits_{\alpha_{1},\alpha_{2}} \dfrac{(\theta_{n}^{(2)} - \theta_{n}^{(1)})^{\alpha_{1}}(v_{1,n}^{(1)} - v_{1,n}^{(2)})^{\alpha_{2}}}{2^{\alpha_{2}}\alpha_{1}!\alpha_{2}!} \nonumber \\
& & +(1- \pi) \sum \limits_{\alpha_{1},\alpha_{2},\beta_{1},\beta_{2}}\dfrac{c^{\alpha_{1}}(c+1)^{\beta_{1}}(\theta_{n}^{(1)} - \theta_{n}^{(2)})^{\alpha_{1}}(\theta_{n}^{(2)})^{\beta_{1}}(v_{2,n}^{(1)} - v_{2,n}^{(2)})^{\alpha_{2}}(v_{2,n}^{(2)}  - v_{1,n}^{(2)})^{\beta_{2}}}{2^{\alpha_{2}+\beta_{2}}\alpha_{1}!\alpha_{2}!\beta_{1}!\beta_{2}!} \nonumber \\
\end{eqnarray}
for all $1 \leq \ell \leq 2\rbar - 2$ where the ranges of $\alpha_{1},\alpha_{2}$ in the first sum of 
$A_{n,\ell}$ satisfy $\alpha_{1}+2\alpha_{2}=\ell$, $1 \leq \alpha_{1}+\alpha_{2} \leq \rbar - 1$ 
and the ranges of $\alpha_{1},\alpha_{2},\beta_{1},\beta_{2}$ in the second sum of $A_{n,\ell}$ satisfy 
$\alpha_{1}+\beta_{1}+2\alpha_{2}+2\beta_{2} = \ell$, $1 \leq \alpha_{1}+\alpha_{2} \leq \rbar - 1$, 
and $0 \leq \beta_{1}+\beta_{2} \leq \rbar - 1 - (\alpha_{1}+\alpha_{2})$.

Now, by our definition of $\bfeta_n^{(1)}, \bfeta_n^{(2)}$ 
and the key equalities \eqref{eq:pf_minimax_asym_identity}, we have
\begin{align*}
A_{n,\ell}= \pi \sum \limits_{\alpha_{1},\alpha_{2}} \dfrac{(\epsilon_n x_1^*)^{\alpha_{1}}(\epsilon_n^2y_1^*)^{\alpha_{2}}}{2^{\alpha_{2}}\alpha_{1}!\alpha_{2}!} 
+(1- \pi) \sum \limits_{\substack{\alpha_{1},\alpha_{2}\\\beta_{1},\beta_{2}}}\dfrac{c^{\alpha_{1}}(c+1)^{\beta_{1}}(-\epsilon_n x_1^*)^{\alpha_{1}}(\epsilon_nx_1^*)^{\beta_{1}}(\epsilon_n^2 y_2^*)^{\alpha_{2}}(\epsilon_n^2 y_3^*)^{\beta_{2}}}{2^{\alpha_{2}+\beta_{2}}\alpha_{1}!\alpha_{2}!\beta_{1}!\beta_{2}!}. 
\end{align*}
Due to the constraint $\alpha_1 + 2\alpha_2 = \ell$ in the first summation, 
and the constraint $\alpha_{1}+\beta_{1}+2\alpha_{2}+2\beta_{2} = \ell$ in the second summation,
the above display reduces to
\begin{align}
\label{eq:pf_minimax_asym_coefficient_reduction}
A_{n,\ell}= \epsilon_n^\ell \left[\pi \sum \limits_{\alpha_{1},\alpha_{2}}
 \dfrac{(x_1^*)^{\alpha_{1}}(y_1^*)^{\alpha_{2}}}{2^{\alpha_{2}}\alpha_{1}!\alpha_{2}!} 
+(1- \pi) \sum \limits_{\substack{\alpha_{1},\alpha_{2}\\\beta_{1},\beta_{2}}}
\dfrac{c^{\alpha_{1}}(c+1)^{\beta_{1}}(-x_1^*)^{\alpha_{1}}(x_1^*)^{\beta_{1}}(y_2^*)^{\alpha_{2}}(y_3^*)^{\beta_{2}}}{2^{\alpha_{2}+\beta_{2}}\alpha_{1}!\alpha_{2}!\beta_{1}!\beta_{2}!}\right]. 
\end{align}
On the other hand, we chose
$(x_1^*, x_2^*, y_1^*, y_2^*, y_3^*)\in \bbR^5$ to be a non-trivial solution to the system of polynomial
equations \eqref{eq:asymmetric_system}.
It follows that $A_{n,\ell} = 0$ for all $\ell = 1, \dots, \rbar-1$.
Also, from \eqref{eq:pf_minimax_asym_coefficient_reduction} we obtain
\begin{align}
\label{eq:pf_minimax_asym_maxcoefs}
\max\Big\{\abss{A_{n,\ell}}: \rbar \leq \ell \leq 2\rbar-2 \Big\} \lesssim \epsilon_n^{\rbar} = \frac 1 n.
\end{align}
We therefore have,\begin{align}
\nonumber
h^2\Big(g( \cdot,\bfeta^{(1)}_n), 
 g(\cdot,\bfeta_n^{(2)})\Big) 
 &= \int \left[\frac{\sum \limits_{\ell=\rbar}^{2\rbar-2} A_{n,\ell} \frac{\partial^\ell f}{\partial\theta^\ell} (x, -\theta_n^{(2)}, v_{1,n}^{(2)}) + R_n(x) }{
 { \sqrt{g(x,\bfeta_n^{(1)})} +
  \sqrt{g(x,\bfeta_n^{(2)})} }}\right]^2 dx \\
 \label{eq:pf_minimax_asym_reduction}
 &\lesssim \int \frac{\sum \limits_{\ell=\rbar}^{2\rbar-2} \left[A_{n,\ell} \frac{\partial^\ell f}{\partial\theta^\ell} (x, -\theta_n^{(2)}, v_{1,n}^{(2)})\right]^2 + R_n^2(x)}{
 {\pi f(x,-\theta_n^{(2)}, v_{1,n}^{(2)}) }} dx.
 \end{align}
It may be verified that for Gaussian densities,
\begin{align}
\label{eq:pf_minimax_asym_integrability}
\int\frac{ \left[\frac{\partial^\ell f}{\partial\theta^\ell} (x, -\theta_n^{(2)}, v_{1,n}^{(2)})\right]^2 }{
 f(x,-\theta_n^{(2)}, v_{1,n}^{(2)})}dx 
  < \infty,\quad \ell=\rbar, \dots, 2\rbar-2.
  \end{align}
Therefore, combining \eqref{eq:pf_minimax_asym_maxcoefs}, \eqref{eq:pf_minimax_asym_reduction} and
\eqref{eq:pf_minimax_asym_integrability}, we obtain
\begin{align}
\label{eq:kl_bound_rem}
h^2\Big(g(&\cdot,\bfeta_n^{(1)}), 
 g(\cdot,\bfeta_n^{(2)})\Big)
  \lesssim \frac  1 n+ \int \frac{R_n^2(x)}{
    {f(x,-\theta_n^{(2)}, v_{1,n}^{(2)}) }}dx. 
\end{align}
Furthermore, notice that
\begin{align}
\label{eq:remainder_decomposition_minimax_asym}
R_n^2(x) \lesssim R_{1,n}^2(x) + R_{2,n}^2(x) +
\sum \limits_{|\alpha| \leq \rbar}  |\theta_{n}^{(1)} - \theta_{n}^{(2)}|^{2\alpha_{1}}|v_{2,n}^{(1)} - v_{2,n}^{(2)}|^{2\alpha_{2}}R_{2,n,\alpha}^2(x).
\end{align}
For $j=1, 2$, we have
\begin{align*}
|R_{j,n}(x)| 
 &\lesssim 
\sum_{|\beta|=\rbar} \frac{|\theta_n^{(2)} - \theta_n^{(1)}|^{\beta_1}|v_{j,n}^{(1)}-v_{j,n}^{(2)}|^{\beta_2}}{\beta_1!\beta_2!} 
 \\ &\times \sup_{0\leq t \leq 1} \left|\frac{\partial^{\rbar} f}{\partial\theta^{\beta_1} \partial v^{\beta_2}}\left(x, -\theta_n^{(2)}+ t(\theta_n^{(1)}-\theta_n^{(2)}), v_{j,n}^{(2)} + t(v_{j,n}^{(1)}-v_{j,n}^{(2)})\right) \right| \\
 &\lesssim \sum_{|\beta|=\rbar}  \epsilon_n^{\beta_1+ \beta_2}
  \sup_{0\leq t \leq 1} \left|\frac{\partial^{\rbar} f}{\partial\theta^{\beta_1} \partial v^{\beta_2}}\left(x, -\theta_n^{(2)}+ t(\theta_n^{(1)}-\theta_n^{(2)}), v_{j,n}^{(2)} + t(v_{j,n}^{(1)}-v_{j,n}^{(2)})\right) \right|,
\end{align*}
and so,
\begin{align}
\nonumber 
 \int & \frac{R_{j,n}^2(x)}{f(x; -\theta_n^{(2)}, v_{j,n}^{(2)})}dx\\
\nonumber 
  &\lesssim \sum_{|\beta|=\rbar} \epsilon_n^{2(\beta_1+ \beta_2)} 
			\int \frac{ 
			  \sup_{0\leq t \leq 1} \left|\frac{\partial^{\rbar} f}{\partial\theta^{\beta_1} \partial v^{\beta_2}}\left(x, -\theta_n^{(2)}+ t(\theta_n^{(1)}-\theta_n^{(2)}), v_{j,n}^{(2)} + t(v_{j,n}^{(1)}-v_{j,n}^{(2)})\right) \right|^2
            }{	
              		{f(x; -\theta_n^{(2)}, v_{j,n}^{(2)})}
            }			dx \\
  &\lesssim \sum_{|\beta|=\rbar} \epsilon_n^{2(\beta_1+ \beta_2)} 
  \lesssim \epsilon_n^{2\rbar} 
  \asymp \frac 1 n.
\label{eq:remainder1_bound}
 \end{align}
By repeating similar calculations for the final term in \eqref{eq:remainder_decomposition_minimax_asym}, 
together with the bound in \eqref{eq:remainder1_bound},
we arrive at
$$\int \frac{R_n^2(x)}{f(x,-\theta_n^{(2)}, v_{1,n}^{(2)}) }dx \lesssim \frac 1 n.$$
Combining the above display with \eqref{eq:kl_bound_rem}, 
then yields
$$h^2\Big(g(\cdot,\bfeta_n^{(1)}), 
 g(\cdot,\bfeta_n^{(2)})\Big) \lesssim \frac 1 n,$$
proving \eqref{eq:hellinger_goal}. The claim follows.  
\hfill $\square$

\subsection{Proof of Proposition \ref{proposition:asymmetric_system}}
\paragraph{PROOF OF PROPOSITION \ref{proposition:asymmetric_system}}
To prove the claim, it suffices to derive a non-trivial solution $(x_1, x_2, y_1, y_2, y_3) \in \bbR^5$
to the system of system of polynomial equations \eqref{eq:asymmetric_system} 
for $r=5$. To this end, we first set $x_1=x_2 = 0$, so that the system reduces to
\begin{align}
\label{eq:pf_prop_asymm_reduction1}
(1-\pi)\sum \limits_{\alpha_{2},\beta_{2}} 
 \dfrac{y_2^{\alpha_2} y_3^{\beta_2}}{2^{\alpha_{2}+\beta_{2}} \alpha_{2}! \beta_{2}!}
+ \pi \sum_{\alpha_2} \dfrac{ y_{1}^{\alpha_2}}{2^{\alpha_2}\alpha_2!} = 0, \quad \ell=1, \dots, 5,
\end{align}
where the first summation in the above display is taken over 
all integers $1 \leq \alpha_2 \leq 5$ and $0 \leq \beta_1 \leq 5 - \alpha_2$
satisfying $2(\alpha_2 + \beta_2) = \ell$, and the second summation is taken
over all integers $1 \leq \alpha_2 \leq 5$ such that $2\alpha_2 = \ell$. Clearly, both of these summations are empty
when $\ell$ is odd, hence equality \eqref{eq:pf_prop_asymm_reduction1} 
holds vacuously for $\ell\in\{1,3,5\}$. It thus remains to show that there exist $y_1, y_2, y_3 \in \bbR$
such that the left-hand side of Eq. \eqref{eq:pf_prop_asymm_reduction1}
vanishes for $\ell\in \{2,4\}$. For such even integers $\ell$, notice that Eq. \eqref{eq:pf_prop_asymm_reduction1}
reduces to
\begin{align*}
(1-\pi)\sum \limits_{\alpha_{2}=1}^{\ell/2}
 \dfrac{y_2^{\alpha_2} y_3^{(\ell/2)-\alpha_2}}{ \alpha_{2}! (\ell/2 - \alpha_2)!}
+ \pi \dfrac{ y_{1}^{\alpha_2}}{(\ell/2)!} = 0, \quad \ell \in \{2,4\},
\end{align*}
which in turn reduces to
\begin{align}
\label{eq:pf_prop_asymm_reduction2}
(1-\pi)\sum \limits_{\alpha_{2}=1}^{\ell/2} {\ell/2 \choose \alpha_2}
 y_2^{\alpha_2} y_3^{\beta_2}
+ \pi y_{1}^{\alpha_2}= 0, \quad \ell \in \{2,4\}.
\end{align}
For $\ell=2$, Eq. \eqref{eq:pf_prop_asymm_reduction2} reads
\begin{align} 
(1-\pi)y_2 + \pi y_1 = 0,
 \end{align}
which is satisfied whenever $y_1 = -y_2/c$. 
Likewise, for $\ell=4$, Eq. \eqref{eq:pf_prop_asymm_reduction2} reads
$$(1-\pi) \big(2y_2y_3 + y_2^2\big) + \pi \left(-\frac{y_2}{c}\right)^2 = 0,$$
which is satisfied whenever $y_3 = -(y_2/2)[1+(1/c)]$, for any $y_2 \in \bbR$. The claim follows. \hfill $\square$

\subsection{Proofs of Lemmas}
\label{sec:lemmas_asym}

\paragraph{PROOF OF LEMMA \ref{lemma:case_a2_asym}} 
We prove the Lemma by considering three cases.
\paragraph{Case a.2.1:} $|v_{2,n}^{(2)} - v_{1,n}^{(2)}|^{1/2}/|\theta_{n}^{(1)} - \theta_{n}^{(2)}| \to \infty$ as $n \to \infty$. From the formulation of $A_{n,3}$, we can easily check that 
\begin{eqnarray}
\dfrac{|A_{n,3}|}{|\theta_{n}^{(1)} - \theta_{n}^{(2)}||v_{2,n}^{(2)} - v_{1,n}^{(2)}|} \to \dfrac{(1-\pi)}{2}. \nonumber
\end{eqnarray}
From the formulation of $D_{n}$, it is clear that
\begin{eqnarray}
\dfrac{D_{n}}{|\theta_{n}^{(1)} - \theta_{n}^{(2)}||v_{2,n}^{(2)} - v_{1,n}^{(2)}|} \to 0. \nonumber
\end{eqnarray}
Therefore, $\max \limits_{1 \leq \ell \leq 2\rbar} |A_{n,\ell}|/D_{n} \not \to 0$. Additionally, 
for each $1 \leq |\alpha| \leq \rbar$, as $n$ is sufficiently large, we have 
\begin{eqnarray}
& & \hspace{- 3 em} \dfrac{|\theta_{n}^{(1)} - \theta_{n}^{(2)}|^{\alpha_{1}}|v_{2,n}^{(1)} - v_{2,n}^{(2)}|^{\alpha_{2}}\|R_{2,n,\alpha}\|_{\infty}}{\max \limits_{1 \leq \ell \leq 2\rbar} |A_{n,\ell}|} \nonumber \\
& & \leq \dfrac{O\parenth{|\theta_{n}^{(1)} - \theta_{n}^{(2)}|^{\alpha_{1}}|v_{2,n}^{(1)} - v_{2,n}^{(2)}|^{\alpha_{2}}(|\theta_{n}^{(2)}|^{\rbar-|\alpha|+\gamma}+|v_{2,n}^{(2)} - v_{1,n}^{(2)}|^{\rbar-|\alpha|+\gamma})}}{|\theta_{n}^{(1)} - \theta_{n}^{(2)}||v_{2,n}^{(2)} - v_{1,n}^{(2)}|}, \nonumber
\end{eqnarray}
which goes to 0 as $n \to \infty$. Hence, we eventually have $\|R_n \|_{\infty}/\max \limits_{1 \leq \ell \leq 2\rbar} |A_{n,\ell}| \to 0$. 
\paragraph{Case a.2.2:} $|v_{2,n}^{(2)} - v_{1,n}^{(2)}|^{1/2}/|\theta_{n}^{(1)} - \theta_{n}^{(2)}| \not \to \infty$ as $n \to \infty$. Equipped with that assumption, we have 
\begin{align}
\max \left\{|v_{1,n}^{(1)} - v_{1,n}^{(2)}|^{1/2}, |v_{2,n}^{(1)} - v_{2,n}^{(2)}|^{1/2}\right\}/|\theta_{n}^{(1)} - \theta_{n}^{(2)}| \to 0. \nonumber
\end{align} 
If $\theta_{n}^{(1)}/\theta_{n}^{(2)} \not \to -1$, then we have $(\theta_{n}^{(1)} + \theta_{n}^{(2)})/(\theta_{n}^{(1)} - \theta_{n}^{(2)}) \not \to 0$. Therefore, we quickly obtain that 
\begin{align}
|A_{n,2}|/|\theta_{n}^{(2)} - \theta_{n}^{(1)}|^{2} \not \to 0. \nonumber
\end{align} 
Since $D_{n}/|\theta_{n}^{(2)} - \theta_{n}^{(1)}|^{2} \to 0$, the previous result implies that $\max \limits_{1 \leq \ell \leq 2\rbar} |A_{n,\ell}|/D_{n} \not \to 0$. Furthermore, for each $1 \leq |\alpha| \leq\rbar$, as $n$ is sufficiently large, we have 
\begin{eqnarray}
& & \hspace{ - 3 em} \dfrac{|\theta_{n}^{(1)} - \theta_{n}^{(2)}|^{\alpha_{1}}|v_{2,n}^{(1)} - v_{2,n}^{(2)}|^{\alpha_{2}}\|R_{2,n,\alpha} \|_{\infty}}{\max \limits_{1 \leq \ell \leq 2\rbar} |A_{n,\ell}|} \nonumber \\
& & \leq \dfrac{O\parenth{|\theta_{n}^{(1)} - \theta_{n}^{(2)}|^{\alpha_{1}}|v_{2,n}^{(1)} - v_{2,n}^{(2)}|^{\alpha_{2}}(|\theta_{n}^{(2)}|^{\rbar-|\alpha|+\gamma}+|v_{2,n}^{(2)} - v_{1,n}^{(2)}|^{\rbar-|\alpha|+\gamma})}}{|\theta_{n}^{(1)} - \theta_{n}^{(2)}|^{2}}, \nonumber
\end{eqnarray}
which goes to 0 for all $1 \leq |\alpha| \leq \rbar$. 
Hence, we have $\|R_n\|_{\infty}/\max \limits_{1 \leq \ell \leq 2\rbar} |A_{n,\ell}| \to 0$.

Overall, we only need to consider the setting that $\theta_{n}^{(1)}/\theta_{n}^{(2)} \to -1$ as $n \to \infty$. Under that setting, we can verify that if $\abss{A_{n,3}}/\abss{\theta_{n}^{(1)}-\theta_{n}^{(2)}}^3 \to 0$, then we have 
\begin{align}
(v_{2,n}^{(2)} - v_{1,n}^{(2)})/\parenth{\theta_{n}^{(1)} - \theta_{n}^{(2)}}^2 \to \dfrac{1-c^2}{8}. \nonumber
\end{align}
However, the above limit leads to
\begin{align}
|A_{n,4}|/|\theta_{n}^{(1)} - \theta_{n}^{(2)}|^{4} \not \to 0. \nonumber
\end{align} 
Therefore, 
$\max\left\{\abss{A_{n,3}},\abss{A_{n,4}}\right\}/\abss{\theta_{n}^{(1)}- \theta_{n}^{(2)}}^4 \not \to 0$ as $n \to \infty$. As $D_{n}/|\theta_{n}^{(1)} - \theta_{n}^{(2)}|^{4} \to 0$, 
the previous result demonstrates that $\max \limits_{1 \leq \ell \leq 2\rbar} |A_{n,\ell}|/D_{n} \not \to 0$. Furthermore, we also have 
\begin{align}
\|R_{2,n,\alpha} \|_{\infty}/\max \limits_{1 \leq \ell \leq 2\rbar} |A_{n,\ell}| \lesssim \|R_{2,n,\alpha} \|_{\infty}/|\theta_{n}^{(1)} - \theta_{n}^{(2)}|^{3} \to 0 \nonumber
\end{align} 
for all $1 \leq |\alpha| \leq \rbar$, which eventually leads to $\|R_n \|_{\infty}/\max \limits_{1 \leq \ell \leq 2\rbar} |A_{n,\ell}| \to 0$. 
\hfill $\square$

\paragraph{PROOF OF LEMMA \ref{lemma:case_a3_asym}} 
To simplify the presentation, we only consider the possibility that 
$\max \left\{|v_{1,n}^{(1)} - v_{1,n}^{(2)}|, |v_{2,n}^{(1)} - v_{2,n}^{(2)}|, |v_{1,n}^{(2)} - v_{2,n}^{(2)}|\right\}/\biggr\{|\theta_{n}^{(1)} - \theta_{n}^{(2)}||\theta_{n}^{(2)}|\biggr\} \not \to \infty$ as $n \to \infty$ since the proof argument for other possibilities 
of this term can be carried out in a similar fashion. 

Under these assumptions, there exist $\bar y_1, \bar y_2 \in \bbR$
such that
$(v_{1,n}^{(1)}-v_{1,n}^{(2)})/\left\{(\theta_{n}^{(2)} - \theta_{n}^{(1)})\theta_{n}^{(2)}\right\} \to \overline{y}_{1}$ and $(v_{2,n}^{(1)} - v_{2,n}^{(2)})/\left\{(\theta_{n}^{(2)} - \theta_{n}^{(1)})\theta_{n}^{(2)}\right\} \to \overline{y}_{2}$ as $n \to \infty$.  We will demonstrate that 
\begin{align}
\max \limits_{1 \leq \ell \leq 2\rbar} |A_{n,\ell}|/\left\{|\theta_{n}^{(1)} - \theta_{n}^{(2)}||\theta_{n}^{(2)}|^{3}\right\} \not \to 0. \nonumber
\end{align} 
Assume by the contrary that $\max \limits_{1 \leq \ell \leq 2\rbar} |A_{n,\ell}|/\left\{|\theta_{n}^{(1)} - \theta_{n}^{(2)}||\theta_{n}^{(2)}|^{3}\right\} \to 0$. By dividing both the numerator and denominator of $|A_{n,\ell}|/\left\{|\theta_{n}^{(1)} - \theta_{n}^{(2)}||\theta_{n}^{(2)}|^{3}\right\}$ by $|\theta_{n}^{(1)} - \theta_{n}^{(2)}||\theta_{n}^{(2)}|^{\ell-1}$ ($2 \leq \ell \leq 4$), as $n \to \infty$, we achieve the following system of polynomial equations
\begin{eqnarray}
\pi \overline{y}_{1} - 2c + (1-\pi)\overline{y}_{2}= 0, \ \overline{y}_{2} = 2c, \ \overline{y}_{2} = 2c(c+1)/3, \nonumber
\end{eqnarray}
which cannot hold as $\pi \in (0,1/2)$. Therefore, 
$\max \limits_{1 \leq \ell \leq 2\rbar} |A_{n,\ell}|/\left\{|\theta_{n}^{(1)} - \theta_{n}^{(2)}||\theta_{n}^{(2)}|^{3}\right\} \not \to 0$. 
From the formulation of $D_{n}$, it is clear that $D_{n}/\left\{|\theta_{n}^{(1)} - \theta_{n}^{(2)}||\theta_{n}^{(2)}|^{3}\right\} \to 0$. 
Hence, $\max \limits_{1 \leq \ell \leq 2\rbar} |A_{n,\ell}|/D_{n} \not \to 0$ as $n \to \infty$. 
Furthermore, we also can check that $\|R_{2,n,\alpha} \|_{\infty}/\left\{|\theta_{n}^{(1)} - \theta_{n}^{(2)}||\theta_{n}^{(2)}|^{3}\right\} \to 0$ 
for all $1 \leq |\alpha| \leq \rbar$. As a consequence, we have $\|R_n \|_{\infty}/\max \limits_{1 \leq \ell \leq \rbar} |A_{n,\ell}| \to 0$. 
The claim follows.
\hfill $\square$ \\

\section{Proofs under the Symmetric Regime}
\label{sec:appendix_symmetric}
\subsection{Proof of Theorem \ref{theorem:convergence_sym} }
Similarly as in the proof of Theorem \ref{theorem:convergence_asym}, 
the key to proving Theorem \ref{theorem:convergence_sym}.(a) is contained in the following result. 
\begin{theorem}
\label{theorem:prelim_sym}
Under the symmetric regime $\pi=1/2$, we have
\begin{eqnarray}
\inf \limits_{\bfeta^{(1)}, \bfeta^{(2)} \in \etaspace} \dfrac{V\parenth{g(\cdot,\bfeta^{(1)}),g(\cdot,\bfeta^{(2)})}}{
\symloss_{\ordersym}^{\ordersym}(\bfeta^{(1)}, \bfeta^{(2)})} > 0. \nonumber
\end{eqnarray}
\end{theorem}

We begin by proving Theorem \ref{theorem:prelim_sym}, and we then prove Theorem \ref{theorem:convergence_sym} below.

\paragraph{PROOF OF THEOREM \ref{theorem:prelim_sym}} 
For simplicity, let $\rbar = \ordersym$ throughout the proof.
Similarly to the proof of Theorem \ref{theorem:prelim_asym}, assume by way of a contradiction that
the claim does not hold. It follows
that we can find sequences $\bfeta_n^{(1)}=(\theta_n^{(1)}, v_{1,n}^{(1)}, v_{2,n}^{(1)}),
\bfeta_n^{(2)}=(\theta_n^{(2)}, v_{1,n}^{(2)}, v_{2,n}^{(2)})$ such that 
\begin{eqnarray}
\label{eq:prelim_sym_assumption}
\frac 1 {\overline D_n} V\parenth{g(\cdot,\bfeta_n^{(1)}), g(\cdot,\bfeta_n^{(2)})} \to 0 ,
\end{eqnarray}
where $\overline{D}_n = \symloss_{\rbar}^{\rbar}(\bfeta_n^{(1)}, \bfeta_n^{(2)})$.
For simplicity of presentation, we 
only consider the most challenging setting where $\theta_{n}^{(i)} \to 0$ ($1 \leq i \leq 2$) 
while $v_{1,n}^{(i)},v_{2,n}^{(i)} \to v_{0}$ ($1 \leq i \leq 2$) for some value $v_{0}>0$. 
Now, we have the following settings with $\theta_{n}^{(1)}, \theta_{n}^{(2)}, v_{1,n}^{(1)}, v_{1,n}^{(2)}, v_{2,n}^{(1)}$, and $v_{2,n}^{(2)}$.
\paragraph{Case b:} $|\theta_{n}^{(1)} - \theta_{n}^{(2)}|^{\rbar}+|v_{1,n}^{(1)} - v_{1,n}^{(2)}|^{\rbar/2} + |v_{2,n}^{(1)} - v_{2,n}^{(2)}|^{\rbar/2} \leq |\theta_{n}^{(1)} + \theta_{n}^{(2)}|^{\rbar}+|v_{1,n}^{(1)} - v_{2,n}^{(2)}|^{\rbar/2} + |v_{1,n}^{(2)} - v_{2,n}^{(1)}|^{\rbar/2}$. Under this setting, we have
\begin{align*}
\overline{D}_{n} = |\theta_{n}^{(1)} - \theta_{n}^{(2)}|^{\rbar}+|v_{1,n}^{(1)} - v_{1,n}^{(2)}|^{\rbar/2} + |v_{2,n}^{(1)} - v_{2,n}^{(2)}|^{\rbar/2}.
\end{align*}
Similarly to the Taylor expansions in the proof of Theorem \ref{theorem:prelim_asym}, 
by Taylor expansion up to order $\rbar$, we obtain
\begin{eqnarray}
& & \hspace{ - 2 em} \frac 1 {\overline D_n} \Big[g(x,\bfeta_n^{(1)}) - g(x,\bfeta_n^{(2)}) \Big] \nonumber \\
& & = \dfrac{1}{2\overline D_n} \left[\parenth{f(x,-\theta_{n}^{(1)},v_{1,n}^{(1)}) - f(x,-\theta_{n}^{(2)},v_{1,n}^{(2)})}+
    \parenth{f(x,\theta_{n}^{(1)},v_{2,n}^{(1)}) - f(x,\theta_{n}^{(2)},v_{2,n}^{(2)})} \right] \nonumber \\
& & = \frac 1 {\overline D_n} \left[\sum \limits_{\ell=1}^{2\rbar} B_{n,\ell}\dfrac{\partial^{\ell}{f}}{\partial{\theta^{\ell}}}(x,-\theta_{n}^{(2)},v_{1,n}^{(2)})+\overline{R}_n(x)\right] \nonumber
\end{eqnarray}
where the formulations of $B_{n,\ell}$ and $\overline{R}_n(x)$ are as follows
\begin{eqnarray}
& & B_{n,\ell}= \dfrac{1}{2} \sum \limits_{\alpha_{1},\alpha_{2}} \dfrac{(\theta_{n}^{(2)} - \theta_{n}^{(1)})^{\alpha_{1}}(v_{1,n}^{(1)} - v_{1,n}^{(2)})^{\alpha_{2}}}{2^{\alpha_{2}}\alpha_{1}!\alpha_{2}!} \nonumber \\
& & +\dfrac{1}{2} \sum \limits_{\alpha_{1},\alpha_{2},\beta_{1},\beta_{2}}\dfrac{2^{\beta_{1}}(\theta_{n}^{(1)} - \theta_{n}^{(2)})^{\alpha_{1}}(\theta_{n}^{(2)})^{\beta_{1}}(v_{2,n}^{(1)} - v_{2,n}^{(2)})^{\alpha_{2}}(v_{2,n}^{(2)}  - v_{1,n}^{(2)})^{\beta_{2}}}{2^{\alpha_{2}+\beta_{2}}\alpha_{1}!\alpha_{2}!\beta_{1}!\beta_{2}!} \nonumber \\
& & \overline{R}_n(x) = \dfrac{1}{2} \overline{R}_{1,n}(x) + \dfrac{1}{2}\overline{R}_{2,n}(x) +
                      \dfrac{1}{2}\sum \limits_{|\alpha| \leq \rbar}\dfrac{1}{2^{\alpha_{2}}} 
                      \dfrac{(\theta_{n}^{(1)} - \theta_{n}^{(2)})^{\alpha_{1}}(v_{2,n}^{(1)} - v_{2,n}^{(2)})^{\alpha_{2}}}{\alpha_{1}!\alpha_{2}!}
                        \overline{R}_{2,n,\alpha}(x) \nonumber
\end{eqnarray}
for all $1 \leq \ell \leq 2\rbar$ where the ranges of $\alpha_{1},\alpha_{2}$ 
in the first sum of $B_{n,\ell}$ satisfy $\alpha_{1}+2\alpha_{2}=\ell$, $1 \leq \alpha_{1}+\alpha_{2} \leq \rbar$ 
and the ranges of $\alpha_{1},\alpha_{2},\beta_{1},\beta_{2}$ in the second sum of $B_{n,\ell}$ 
satisfy $\alpha_{1}+\beta_{1}+2\alpha_{2}+2\beta_{2} = \ell$, $1 \leq \alpha_{1}+\alpha_{2} \leq \rbar$, 
and $0 \leq \beta_{1}+\beta_{2} \leq \rbar - (\alpha_{1}+\alpha_{2})$. 
Additionally, $\overline{R}_{1,n}(x)$ is Taylor remainder from expanding 
$f(x,-\theta_{n}^{(1)},v_{1,n}^{(1)})$ around $f(x,-\theta_{n}^{(2)},v_{1,n}^{(2)})$ 
up to the $\rbar$ order, $\overline{R}_{2,n}(x)$ is Taylor remainder from expanding 
$f(x,\theta_{n}^{(1)},v_{1,n}^{(2)})$ around $f(x,\theta_{n}^{(2)},v_{2,n}^{(2)})$ up 
to the $\rbar$ order, and $\overline{R}_{2,n,\alpha}(x)$ is Taylor remainder from expanding 
$\dfrac{\partial^{\alpha_{1}+2\alpha_{2}}{f}}{\partial{\theta^{\alpha_{1}+2\alpha_{2}}}}(x,\theta_{n}^{(2)},v_{2,n}^{(2)})$ 
around $\dfrac{\partial^{\alpha_{1}+2\alpha_{2}}{f}}{\partial{\theta^{\alpha_{1}+2\alpha_{2}}}}(x;-\theta_{n}^{(2)},v_{1,n}^{(2)})$ 
up to the order $\rbar-|\alpha|$ for $1 \leq |\alpha| \leq \rbar$. 

Similarly to the proof of Theorem \ref{theorem:prelim_asym}, we have the following settings under Case b.
\paragraph{Case b.1:} $\theta_{n}^{(2)}/\theta_{n}^{(1)} \not \to 1$ and $|v_{2,n}^{(2)} - v_{1,n}^{(2)}|/\max \left\{|v_{1,n}^{(1)} - v_{1,n}^{(2)}|, |v_{2,n}^{(1)} - v_{2,n}^{(2)}|\right\} \not \to \infty$ as $n \to \infty$. 

Assume by way of a contradiction  that all the coefficients $B_{n,\ell}/\overline{D}_{n} \to 0$ for all $1 \leq \ell \leq 2\rbar$. 
We again write
\begin{eqnarray}
\overline{M}_{n} = \left\{|\theta_{n}^{(2)} - \theta_{n}^{(1)}|, |v_{1,n}^{(1)} - v_{1,n}^{(2)}|^{1/2}, |v_{2,n}^{(1)} - v_{2,n}^{(2)}|^{1/2}\right\}. \nonumber
\end{eqnarray}
From the assumption of Case b.1, we have $|v_{2,n}^{(2)} - v_{1,n}^{(2)}|/\overline{M}_{n}^{2} \not \to \infty$ and $|\theta_{n}^{(2)}|/\overline{M}_{n} \not \to \infty$. Therefore, we can define 
$$(\theta_{n}^{(2)} - \theta_{n}^{(1)})/\overline{M}_{n} \to x_{1}, \quad \theta_{n}^{(2)}/\overline{M}_{n} \to x_{2},$$
and,
$$(v_{1,n}^{(1)} - v_{1,n}^{(2)})/\overline{M}_{n}^{2} \to y_{1}, \quad 
  (v_{2,n}^{(1)} - v_{2,n}^{(2)})/\overline{M}_{n}^{2} \to y_{2}, \quad
  (v_{2,n}^{(2)} - v_{1,n}^{(2)})/\overline{M}_{n}^{2} \to y_{3}$$
From the definition of $\overline{M}_{n}$, at least one of $x_{1},y_{1},y_{2}$ is different from 0. 
Additionally, the definition of $x_{1},x_{2}, y_{1},y_{2}$, and $y_{3}$ 
leads to $(\theta_{n}^{(1)}+ \theta_{n}^{(2)})/\overline{M}_{n} \to 2x_{2}-x_{1}$, $(v_{2,n}^{(2)}-v_{1,n}^{(1)})/\overline{M}_{n}^2 \to y_{3}-y_{1}$, 
and $(v_{2,n}^{(1)} - v_{1,n}^{(2)})/\overline{M}_{n}^2 \to y_{2}+y_{3}$. 
By dividing both sides of the assumption of Case b.1
assumption by $\overline{M}_{n}^{\rbar}$ and let $n \to \infty$, we obtain the following constraint with $x_{1},x_{2},y_{1},y_{2},y_{3}$
\begin{eqnarray}
\hspace{ 2 em} \abss{x_{1}}^{\rbar}+\abss{y_{1}}^{\rbar/2}+\abss{y_{2}}^{\rbar/2} \leq \abss{2x_{2}-x_{1}}^{\rbar}+\abss{y_{3}-y_{1}}^{\rbar/2}+\abss{y_{2}+y_{3}}^{\rbar/2}. \label{eqn:system_inequality_symmetric_case}
\end{eqnarray} 
Now, by dividing both the numerator and the denominator of $B_{n,\ell}$ ($1 \leq \ell \leq \rbar$) by $\overline{M}_{n}^{\ell}$, as $n \to \infty$, we have the following system of polynomial equations
\begin{eqnarray}
\hspace{ - 3 em} \sum \limits_{\alpha_{1},\alpha_{2},\beta_{1},\beta_{2}} \dfrac{1}{2^{\alpha_{2}+\beta_{2}}} \dfrac{2^{\beta_{1}}(-x_{1})^{\alpha_{1}}x_{2}^{\beta_{1}}y_{2}^{\alpha_{2}}y_{3}^{\beta_{2}}}{\alpha_{1}!\alpha_{2}!\beta_{1}!\beta_{2}!} + \sum \limits_{\alpha_{1},\alpha_{2}} \dfrac{1}{2^{\alpha_{2}}}\dfrac{x_{1}^{\alpha_{1}}y_{1}^{\alpha_{2}}}{\alpha_{1}!\alpha_{2}!}  = 0 \label{eqn:system_polynomial_symmetric_case}
\end{eqnarray} 
as $\ell=1,\ldots,\rbar$. The above system of polynomial equations along 
with inequality \eqref{eqn:system_inequality_symmetric_case} forms a 
semialgebraic set with the constraint that at least one of $x_{1},y_{1},y_{2}$ 
is different from 0. According to the definition of $\rbar$, 
this semialgebraic set is empty, which is a contradiction. 
Therefore, not all the coefficients $B_{n,\ell}/\overline{D}_{n}$ 
go to 0 as $n \to \infty$. Denote 
$m_{n} = \overline{D}_{n}/\max \limits_{1 \leq \ell \leq 2\rbar} |B_{n,\ell}|$. 
Governed by the previous result, we have $m_{n} \not \to \infty$. Now, we have that
\begin{eqnarray}
m_{n} \biggr(\sum \limits_{\ell=1}^{2\rbar} B_{n,\ell}\dfrac{\partial^{\ell}{f}}{\partial{\theta^{\ell}}}(x,-\theta_{n}^{(2)},v_{1,n}^{(2)})\biggr)/D_{n} \to \sum \limits_{\ell=1}^{2\rbar} \bar{\tau}_{\ell}\dfrac{\partial^{\ell}{f}}{\partial{\theta^{\ell}}}(x,0,v_{0}) \nonumber
\end{eqnarray}
for some coefficients $\bar{\tau}_{\ell}$ which are not all zero. 
Similarly as in the proof of Theorem \ref{theorem:prelim_asym},
by means of Fatou's lemma combined with the hypothesis \eqref{eq:prelim_sym_assumption}, the following then holds
\begin{align}
\sum \limits_{\ell=1}^{2\rbar} \bar{\tau}_{\ell}\dfrac{\partial^{\ell}{f}}{\partial{\theta^{\ell}}}(x,0,v_{0}) = 0. \nonumber
\end{align} 
However, Lemma \ref{lemma:gaussian_si} implies $\bar{\tau}_{\ell} = 0$ for all $1 \leq \ell \leq 2\rbar$, which is a contradiction. 
Therefore, Case b.1 does not hold. 

\paragraph{Case b.2:} $\theta_{n}^{(2)}/\theta_{n}^{(2)} \not \to 1$ and $|v_{2,n}^{(2)} - v_{1,n}^{(2)}|/\max \left\{|v_{1,n}^{(1)} - v_{1,n}^{(2)}|, |v_{2,n}^{(1)} - v_{2,n}^{(2)}|\right\} \to \infty$ as $n \to \infty$. Following the strategy of Case a.2 in part (a), 
it follows from the following Lemma that Case b.2 cannot hold.
\begin{lemma}
\label{lemma:case_b2_sym}
Under the setting of Case b.2, we have
$$\max \limits_{1 \leq \ell \leq 2\rbar} |B_{n,\ell}|/\overline{D}_{n} \not \to 0, \quad\text{and,}\quad
\|\overline{R}_n\|_{\infty}/\max \limits_{1 \leq \ell \leq 2\rbar} |B_{n,\ell}| \to 0.$$
\end{lemma}
The proof of Lemma \ref{lemma:case_b2_sym} appears in Appendix \ref{sec:proofs_lemmas_sym}. 

\paragraph{Case b.3:} $\theta_{n}^{(2)}/\theta_{n}^{(1)} \to 1$ as $n \to \infty$. 
Once again, it follows from the following Lemma that Case b.3 cannot hold.
\begin{lemma}
\label{lemma:case_b3_sym}
Under the setting of Case b.3, we have
$$\max \limits_{1 \leq \ell \leq 2\rbar} |B_{n,\ell}|/\overline{D}_{n} \not \to 0,\quad
\text{and,}\quad\|\overline{R}_n\|_{\infty}/\max \limits_{1 \leq \ell \leq 2\rbar} |B_{n,\ell}| \to 0.$$
\end{lemma}
The proof of Lemma \ref{lemma:case_b3_sym} appears in Appendix \ref{sec:proofs_lemmas_sym}. 
Altogether, we conclude that case $b$ cannot hold.

\paragraph{Case c:} $|\theta_{n}^{(1)} - \theta_{n}^{(2)}|^{\rbar}+|v_{1,n}^{(1)} - v_{1,n}^{(2)}|^{\rbar/2} + |v_{2,n}^{(1)} - v_{2,n}^{(2)}|^{\rbar/2} > |\theta_{n}^{(1)} + \theta_{n}^{(2)}|^{\rbar}+|v_{1,n}^{(1)} - v_{2,n}^{(2)}|^{\rbar/2} + |v_{1,n}^{(2)} - v_{2,n}^{(1)}|^{\rbar/2}$. 

Under this setting, we have
\begin{eqnarray}
\overline{D}_{n} = |\theta_{n}^{(1)} + \theta_{n}^{(2)}|^{\rbar}+|v_{1,n}^{(1)} - v_{2,n}^{(2)}|^{\rbar/2} + |v_{1,n}^{(2)} - v_{2,n}^{(1)}|^{\rbar/2}. \nonumber
\end{eqnarray}
Similarly to Case b, by means of Taylor expansion up to the $\rbar$ order, we obtain that
\begin{align*}
\frac 1 {\overline{D}_n}&\Big[g(x,\bfeta_n^{(1)}) - g(x,\bfeta_n^{(2)})\Big] \\
 &= \frac{1}{2\overline D_n} 
   \Big[\Big(f(x,-\theta_{n}^{(1)},v_{1,n}^{(1)}) - f(x,\theta_{n}^{(2)},v_{2,n}^{(2)})\Big)+
        \Big(f(x,\theta_{n}^{(1)},v_{2,n}^{(1)}) - f(x,-\theta_{n}^{(2)},v_{1,n}^{(2)})\Big)\Big]  \\
 &= \frac 1 {\overline D_n} \left[\sum_{\ell=1}^{2\rbar} C_{n,\ell}\dfrac{\partial^{\ell}{f}}{\partial{\theta^{\ell}}}(x,-\theta_{n}^{(2)},v_{1,n}^{(2)})+\widetilde{R}_n(x)\right]
\end{align*}
where,
\begin{eqnarray}
& & C_{n,\ell}= \dfrac{1}{2} \sum \limits_{\alpha_{1},\alpha_{2}} \dfrac{(-\theta_{n}^{(2)} - \theta_{n}^{(1)})^{\alpha_{1}}(v_{1,n}^{(1)} - v_{2,n}^{(2)})^{\alpha_{2}}}{2^{\alpha_{2}}\alpha_{1}!\alpha_{2}!} \nonumber \\
& & +\dfrac{1}{2} \sum \limits_{\alpha_{1},\alpha_{2},\beta_{1},\beta_{2}}\dfrac{2^{\beta_{1}}(\theta_{n}^{(1)} + \theta_{n}^{(2)})^{\alpha_{1}}(-\theta_{n}^{(2)})^{\beta_{1}}(v_{2,n}^{(1)} - v_{1,n}^{(2)})^{\alpha_{2}}(v_{1,n}^{(2)}  - v_{2,n}^{(2)})^{\beta_{2}}}{2^{\alpha_{2}+\beta_{2}}\alpha_{1}!\alpha_{2}!\beta_{1}!\beta_{2}!} \nonumber \\
& & \widetilde{R}_n(x) = \dfrac{1}{2} \widetilde{R}_{1,n}(x) + \dfrac{1}{2}\widetilde{R}_{2,n}(x) + \dfrac{1}{2}\sum \limits_{|\alpha| \leq \rbar}\dfrac{1}{2^{\alpha_{2}}} \dfrac{(\theta_{n}^{(1)} + \theta_{n}^{(2)})^{\alpha_{1}}(v_{2,n}^{(1)} - v_{1,n}^{(2)})^{\alpha_{2}}}{\alpha_{1}!\alpha_{2}!}\widetilde{R}_{2,n,\alpha}(x) \nonumber
\end{eqnarray}
for all $1 \leq \ell \leq 2\rbar$ where the ranges of $\alpha_{1},\alpha_{2}$ 
in the first sum of $C_{n,\ell}$ satisfy 
$\alpha_{1}+2\alpha_{2}=\ell$, $1 \leq \alpha_{1}+\alpha_{2} \leq \rbar$ 
and the ranges of $\alpha_{1},\alpha_{2},\beta_{1},\beta_{2}$ 
in the second sum of $C_{n,\ell}$ satisfy $\alpha_{1}+\beta_{1}+2\alpha_{2}+2\beta_{2} = \ell$, 
$1 \leq \alpha_{1}+\alpha_{2} \leq \rbar$, and $0 \leq \beta_{1}+\beta_{2} \leq \rbar - (\alpha_{1}+\alpha_{2})$. 
Additionally, $\widetilde{R}_{1,n}(x)$ is Taylor remainder from expanding $f(x,-\theta_{n}^{(1)},v_{1,n}^{(1)})$ 
around $f(x,\theta_{n}^{(2)},v_{2,n}^{(2)})$ up to the $\rbar$ order, $\widetilde{R}_{2,n}(x)$ 
is Taylor remainder from expanding $f(x,\theta_{n}^{(1)},v_{1,n}^{(2)})$ around $f(x,-\theta_{n}^{(2)},v_{1,n}^{(2)})$ 
up to the $\rbar$ order, and $\widetilde{R}_{2,n,\alpha}(x)$ is Taylor remainder from 
expanding $\dfrac{\partial^{\alpha_{1}+2\alpha_{2}}{f}}{\partial{\theta^{\alpha_{1}+2\alpha_{2}}}}(x,-\theta_{n}^{(2)},v_{1,n}^{(2)})$ 
around $\dfrac{\partial^{\alpha_{1}+2\alpha_{2}}{f}}{\partial{\theta^{\alpha_{1}+2\alpha_{2}}}}(x,\theta_{n}^{(2)},v_{2,n}^{(2)})$ up 
to the order $\rbar-|\alpha|$ for $1 \leq |\alpha| \leq \rbar$. 

To ease the proof argument, we only consider the setting that $\theta_{n}^{(2)}/\theta_{n}^{(1)} \not \to 1$ and $|v_{2,n}^{(2)} - v_{1,n}^{(2)}|/\max \left\{|v_{1,n}^{(1)} - v_{2,n}^{(2)}|, |v_{2,n}^{(1)} - v_{1,n}^{(2)}|\right\} \not \to \infty$ as $n \to \infty$. The other possibilities of these terms can be argued similarly as those in Case b.2 and Case b.3.  Assume now that all the coefficients $C_{n,\ell}/\overline{D}_{n} \to 0$ for all 
$1 \leq \ell \leq 2\rbar$. Denote
\begin{eqnarray}
\widetilde{M}_{n} = \left\{|\theta_{n}^{(2)} - \theta_{n}^{(1)}|, |v_{1,n}^{(1)} - v_{2,n}^{(2)}|^{1/2}, |v_{2,n}^{(1)} - v_{1,n}^{(2)}|^{1/2}\right\}. \nonumber
\end{eqnarray}
From the previous assumptions, we have $|v_{2,n}^{(2)} - v_{1,n}^{(2)}|/\widetilde{M}_{n}^{2} \not \to \infty$ and 
$|\theta_{n}^{(2)}|/\widetilde{M}_{n} \not \to \infty$. 
Therefore, we define 
$$(-\theta_{n}^{(2)} - \theta_{n}^{(1)})/\widetilde{M}_{n} \to \bar{x}_{1}, \quad
  -\theta_{n}^{(2)}/\widetilde{M}_{n} \to \bar{x}_{2}, $$
$$  (v_{1,n}^{(1)} - v_{2,n}^{(2)})/\widetilde{M}_{n}^{2} \to \bar{y}_{1},\quad
  (v_{2,n}^{(1)} - v_{1,n}^{(2)})/\widetilde{M}_{n}^{2} \to \bar{y}_{2}, \quad
  (v_{1,n}^{(2)} - v_{2,n}^{(2)})/\widetilde{M}_{n}^{2} \to \bar{y}_{3}.$$
According to the definition of $\widetilde{M}_{n}$, at least one of $\bar{x}_{1},\bar{y}_{1},\bar{y}_{2}$ is different from 0. Additionally, the definition of $\bar{x}_{1},\bar{x}_{2}, \bar{y}_{1},\bar{y}_{2}$, and $\bar{y}_{3}$ leads to $(\theta_{n}^{(1)}- \theta_{n}^{(2)})/\widetilde{M}_{n} \to 2\bar{x}_{2}-\bar{x}_{1}$, $(v_{2,n}^{(1}-v_{2,n}^{(2)})/\widetilde{M}_{n}^2 \to \bar{y}_{2}+\bar{y}_{3}$, and $(v_{1,n}^{(1)} - v_{1,n}^{(2)})/\widetilde{M}_{n}^2 \to \bar{y}_{1}-\bar{y}_{3}$. According to the assumption of Case c, by dividing both sides of this assumption by $\widetilde{M}_{n}^{\rbar}$ and let $n \to \infty$, we obtain the following inequality
\begin{eqnarray}
\hspace{4 em} \abss{\bar{x}_{1}}^{\rbar}+\abss{\bar{y}_{1}}^{\rbar/2}+\abss{\bar{y}_{2}}^{\rbar/2} \leq \abss{2\bar{x}_{2}-\bar{x}_{1}}^{\rbar}+\abss{\bar{y}_{3}-\bar{y}_{1}}^{\rbar/2}+\abss{\bar{y}_{2}+\bar{y}_{3}}^{\rbar/2}. \label{eqn:system_inequality_symmetric_case_second}
\end{eqnarray} 
Now, by dividing both the numerator and the denominator of $C_{n,\ell}$ ($1 \leq \ell \leq \rbar$) by $\widetilde{M}_{n}^{\ell}$, as $n \to \infty$, we have the following system of polynomial equations
\begin{eqnarray}
\hspace{ - 3 em} \sum \limits_{\alpha_{1},\alpha_{2},\beta_{1},\beta_{2}} \dfrac{1}{2^{\alpha_{2}+\beta_{2}}} \dfrac{2^{\beta_{1}}(-\bar{x}_{1})^{\alpha_{1}}\bar{x}_{2}^{\beta_{1}}\bar{y}_{2}^{\alpha_{2}}\bar{y}_{3}^{\beta_{2}}}{\alpha_{1}!\alpha_{2}!\beta_{1}!\beta_{2}!} + \sum \limits_{\alpha_{1},\alpha_{2}} \dfrac{1}{2^{\alpha_{2}}}\dfrac{\bar{x}_{1}^{\alpha_{1}}\bar{y}_{1}^{\alpha_{2}}}{\alpha_{1}!\alpha_{2}!}  = 0 \label{eqn:system_polynomial_symmetric_case_second}
\end{eqnarray}
for all $\ell=1,\ldots,\rbar$. According to the definition of $\rbar$, the system of polynomial equations \eqref{eqn:system_polynomial_symmetric_case_second} and inequality \eqref{eqn:system_inequality_symmetric_case_second} cannot hold unless $\bar{x}_{1}=\bar{y}_{1}=\bar{y}_{2}=0$, which is a contradiction. 
Therefore, not all the coefficients $C_{n,\ell}/\overline{D}_{n}$ go to 0 as $n \to \infty$. 
As a consequence, by means of Fatou's argument, we deduce that Case c cannot happen. The claim follows.
\hfill$\square$

We are now in a position to prove Theorem \ref{theorem:convergence_sym}.

\paragraph{PROOF OF THEOREM \ref{theorem:convergence_sym}}
We fix $\rbar = \ordersym$ throughout the proof.

(a) Similarly as in the proof of Theorem \ref{theorem:convergence_asym}.(a), 
Theorem \ref{theorem:prelim_sym} implies the existence of a universal constant $C_1 > 0$, depending only on $\Omega, \Theta$, such that
for all $\bfeta \in \etaspace$,
\begin{align*}
\symloss_{\rbar}^{\rbar}(\hbfeta_n, \bfeta)
   \leq C_1 h\big(g(\cdot,\hbfeta_n),g(\cdot,\bfeta)\big),
\end{align*}
The claim then follows by an application of Lemma \ref{lemma:hellinger_rate}.
 
(b)  The proof follows along similar lines as that of Theorem \ref{theorem:convergence_asym}.(b).
By definition of $\rbar$, there exists a solution
$(x_1^*, x_2^*, y_1^*, y_2^*, y_3^*) \in \bbR^5$  to the system
of polynomial equalities and inequalities \eqref{eq:symmetric_system_equalities} and \eqref{eq:symmetric_system_inequalities}, with respect to the choice $r = \rbar-1$. 
Set $\bfeta_n^{(i)} = (\theta_n^{(i)}, v_{1,n}^{(i)}, v_{2,n}^{(i)}) \in H$
for $i=1,2$, where, for $\epsilon_n = n^{-1/2\rbar}$,
$$\theta_n^{(1)} = \epsilon_n (x_2^*-x_1^*),~~
  \theta_n^{(2)} = \epsilon_n x_2^*,$$
and where, 
$$ v_{1,n}^{(1)} = \epsilon_n^2 y_1^* + v_0 , ~~
   v_{2,n}^{(1)} = \epsilon_n^2 (y_2^* + y_3^*) + v_0 , ~~
   v_{1,n}^{(2)} = v_0 , ~~
   v_{2,n}^{(2)} = \epsilon_n^2 y_3^* + v_0.$$

The definition of $\epsilon_n$ then implies
that $\assymloss_{\rbar}(\bfeta_n^{(j)},\bfeta_0) \leq c_2 n^{-1/2\rbar}$ for some $c_2 > 0$, for $j=1,2$. 
Furthermore, this choice of parameters satisfies the identities
\begin{equation}
\label{eq:pf_minimax_sym_identity} 
\theta_n^{(2)} - \theta_n^{(1)} = \epsilon_n x_1^*, \quad
  v_{1,n}^{(1)} - v_{1,n}^{(2)} = \epsilon_n^2 y_1^*, \quad
  v_{2,n}^{(1)} - v_{2,n}^{(2)} = \epsilon_n^2 y_2^*,\quad
  v_{2,n}^{(2)} - v_{1,n}^{(2)} = \epsilon_n^2 y_3^*.
  \end{equation}
Invoking Le Cam's Inequality (\cite{tsybakov2008}, Theorem 2.2), we obtain
\begin{align}
\label{eq:pf_sym_minimax_le_cam}
\inf_{\hat \bfeta_n} &\sup_{\bfeta \in \etaspace(c_2n^{-1/2})}
 \bbE_{\bfeta} \big[\assymloss_{\rbar}(\hbfeta_n, \bfeta) \big]
 \geq \assymloss_{\rbar}(\bfeta_n^{(1)}, \bfeta_n^{(2)}) \Big(1 - V\big(g(\cdot, \bfeta^{(1)}_n), g(\cdot, \bfeta_n^{(2)})\big)\Big)
\end{align}
Notice that the identities \eqref{eq:pf_minimax_sym_identity} imply
\begin{align}
\label{eq:pf_minimax_sym_minimum}
\assymloss_{\rbar}^{\rbar}(\bfeta_n^{(1)}, \bfeta_n^{(2)}) 
 &= \epsilon_n^{\rbar} \min \Big\{ |x_1^*|^{\rbar} + |y_1^*|^{\rbar/2} + |y_2^*|^{\rbar/2}, 
 							 |2x_2^* - x_1|^{\rbar} + |y_3^*-y_1^*|^{\rbar/2} + |y_2^* + y_3^*|^{\rbar /2}  \Big\}.
\end{align}
We will argue that the minimum in the above display is nonzero. 
To this end, since $(x_1^*$, $x_2^*,$ $y_1^*, y_2^*$, $y_3^*)$ form a non-trivial solution to the system
of polynomial equations and inequalities in \eqref{eq:symmetric_system_equalities} and
\eqref{eq:symmetric_system_inequalities}, it must hold that one of $x_1^*, y_1^*, y_2^*$
is nonzero, and in particular,
\begin{align*}
0 < |x_1^*|^{\rbar-1} + |y_1^*|^{\frac{\rbar-1}{2}} + |y_2^*|^{\frac{\rbar-1}{2}}.
\end{align*}
 This fact combined with inequality \eqref{eq:symmetric_system_inequalities}
 of the asymmetric system
implies
\begin{align*}
0 < 	 |2x_2^* - x_1^*|^{\rbar-1} + |y_3^*-y_1^*|^{\frac{\rbar-1}{2}} + |y_2^* + y_3^*|^{\frac{\rbar -1}{2}},
\end{align*} 					
from which it follows that the minimum in equation \eqref{eq:pf_minimax_sym_minimum}
is strictly positive. Therefore, we have $\assymloss_{\rbar}(\bfeta_n^{(1)}, \bfeta_n^{(2)}) \asymp \epsilon_n=n^{-1/2\rbar}$. 
Returning to equation \eqref{eq:pf_sym_minimax_le_cam}, and 
using the inequality $V \leq h$ together with the tensorization property of the Hellinger distance (\cite{tsybakov2008}, p. 83),
we obtain
\begin{align*}
\inf_{\hat \bfeta_n} \sup_{\bfeta \in \etaspace(c_2n^{-1/2})}& \bbE_{\bfeta}
 \big[\assymloss_{\rbar}(\hbfeta_n, \bfeta) \big]\\
 &\gtrsim  n^{-1/{2\rbar}}
 \left( 1 - \sqrt{1-\left[1-h^2\Big(g(\cdot,\bfeta_n^{(1)}), g(\cdot,\bfeta_n^{(2)})\Big)\right]^n}\right).
 \end{align*}
Notice that the right-hand side of the above display will be of order $n^{-1/2 }$ provided 
\begin{equation}
\label{eq:pf_minimax_sym_hellinger_goal}
h^2\Big(g(\cdot,\bfeta_n^{(1)}), g(\cdot,\bfeta_n^{(2)})\Big) \lesssim \frac 1 n.
\end{equation}
To prove the claim, it will therefore suffice to prove that \eqref{eq:pf_minimax_sym_hellinger_goal} holds. 
We argue similarly as in the proof of Theorem \ref{theorem:convergence_asym}(b).
Notice that
\begin{equation}
h^2\Big(g(\cdot,\bfeta_n^{(1)}), 
 g(\cdot,\bfeta_n^{(2)})\Big)
 = \int \frac{\left[g(x,\bfeta_n^{(1)}) -
 g(\cdot, \bfeta_n^{(2)})\right]^2}{
 {\left[\sqrt{g(x,\bfeta_n^{(1)})} +
  \sqrt{g(x,\bfeta_n^{(2)})}\right]^2 }} dx
 \end{equation}
We begin by analyzing the numerator of the integrand in the above display.
By a similar Taylor expansion as in Case b of Theorem \ref{theorem:prelim_sym},
but now up to order $\rbar-1$,
we have
\begin{align*}
g(x,&\bfeta_n^{(1)}) - g(x,\bfeta_n^{(2)}) \\
 &= \frac 1 2 \parenth{f(x,-\theta_{n}^{(1)},v_{1,n}^{(1)}) - f(x,-\theta_{n}^{(2)},v_{1,n}^{(2)})}+
    \frac 1 2 \parenth{f(x,\theta_{n}^{(1)},v_{2,n}^{(1)}) - f(x,\theta_{n}^{(2)},v_{2,n}^{(2)})}  \\
 &=\sum \limits_{\ell=1}^{2(\rbar-1)} B_{n,\ell}\dfrac{\partial^{\ell}{f}}{\partial{\theta^{\ell}}}(x,-\theta_{n}^{(2)},v_{1,n}^{(2)})+R_n(x),
\end{align*}
where 
\begin{align*}
B_{n,\ell}&= \dfrac{1}{2} \sum \limits_{\alpha_{1},\alpha_{2}} \dfrac{(\theta_{n}^{(2)} - \theta_{n}^{(1)})^{\alpha_{1}}(v_{1,n}^{(1)} - v_{1,n}^{(2)})^{\alpha_{2}}}{2^{\alpha_{2}}\alpha_{1}!\alpha_{2}!} \\
& +\dfrac{1}{2} \sum \limits_{\alpha_{1},\alpha_{2},\beta_{1},\beta_{2}}\dfrac{2^{\beta_{1}}(\theta_{n}^{(1)} - \theta_{n}^{(2)})^{\alpha_{1}}(\theta_{n}^{(2)})^{\beta_{1}}(v_{2,n}^{(1)} - v_{2,n}^{(2)})^{\alpha_{2}}(v_{2,n}^{(2)}  - v_{1,n}^{(2)})^{\beta_{2}}}{2^{\alpha_{2}+\beta_{2}}\alpha_{1}!\alpha_{2}!\beta_{1}!\beta_{2}!}   \\
 R_n(x) &= \dfrac{1}{2} R_{1,n}(x) + \dfrac{1}{2}R_{2,n}(x) +
                      \dfrac{1}{2}\sum \limits_{|\alpha| \leq \rbar-1}\dfrac{1}{2^{\alpha_{2}}} 
                      \dfrac{(\theta_{n}^{(1)} - \theta_{n}^{(2)})^{\alpha_{1}}(v_{2,n}^{(1)} - v_{2,n}^{(2)})^{\alpha_{2}}}{\alpha_{1}!\alpha_{2}!}
                        R_{2,n,\alpha}(x),
\end{align*}
for all $1 \leq \ell \leq 2(\rbar-1)$ where the ranges of $\alpha_{1},\alpha_{2}$ 
in the first sum of $B_{n,\ell}$ satisfy $\alpha_{1}+2\alpha_{2}=\ell$, $1 \leq \alpha_{1}+\alpha_{2} \leq \rbar-1$ 
and the ranges of $\alpha_{1},\alpha_{2},\beta_{1},\beta_{2}$ in the second sum of $B_{n,\ell}$ 
satisfy $\alpha_{1}+\beta_{1}+2\alpha_{2}+2\beta_{2} = \ell$, $1 \leq \alpha_{1}+\alpha_{2} \leq \rbar-1$, 
and $0 \leq \beta_{1}+\beta_{2} \leq \rbar - 1 - (\alpha_{1}+\alpha_{2})$. 
Further, for $j=1,2$, $R_{j,n}(x)$ is the Taylor remainder arising from an expansion of  
$f(x,-\theta_{n}^{(1)},v_{j,n}^{(1)})$ around $f(x,-\theta_{n}^{(2)},v_{j,n}^{(2)})$ 
up to order $\rbar$,  and $R_{2,n,\alpha}(x)$ is Taylor remainder arising from an expansion of 
$\dfrac{\partial^{\alpha_{1}+2\alpha_{2}}{f}}{\partial{\theta^{\alpha_{1}+2\alpha_{2}}}}(x,\theta_{n}^{(2)},v_{2,n}^{(2)})$ 
around $\dfrac{\partial^{\alpha_{1}+2\alpha_{2}}{f}}{\partial{\theta^{\alpha_{1}+2\alpha_{2}}}}(x;-\theta_{n}^{(2)},v_{1,n}^{(2)})$ 
up to order $\rbar-1-|\alpha|$, for $1 \leq |\alpha| \leq \rbar-1$. 
	
Now, since $\bfeta_n^{(1)}, \bfeta_n^{(2)}$ satisfy the identities \eqref{eq:pf_minimax_sym_identity}, we have
\begin{align*}
B_{n,\ell}
 &=  \sum \limits_{\alpha_{1},\alpha_{2}} \dfrac{(\epsilon_n x_1^*)^{\alpha_{1}}(\epsilon_n^2y_1^*)^{\alpha_{2}}}{2^{\alpha_{2}}\alpha_{1}!\alpha_{2}!} 
+ \sum \limits_{\substack{\alpha_{1},\alpha_{2}\\\beta_{1},\beta_{2}}}\dfrac{ 2^{\beta_{1}}(-\epsilon_n x_1^*)^{\alpha_{1}}(\epsilon_nx_2^*)^{\beta_{1}}(\epsilon_n^2 y_2^*)^{\alpha_{2}}(\epsilon_n^2 y_3^*)^{\beta_{2}}}{2^{\alpha_{2}+\beta_{2}}\alpha_{1}!\alpha_{2}!\beta_{1}!\beta_{2}!} \\
 &= \epsilon_n^\ell \left[\sum \limits_{\alpha_{1},\alpha_{2}} \dfrac{(x_1^*)^{\alpha_{1}}( y_1^*)^{\alpha_{2}}}{2^{\alpha_{2}}\alpha_{1}!\alpha_{2}!} 
+ \sum \limits_{\substack{\alpha_{1},\alpha_{2}\\\beta_{1},\beta_{2}}}\dfrac{ 2^{\beta_{1}}(- x_1^*)^{\alpha_{1}}( x_2^*)^{\beta_{1}}(y_2^*)^{\alpha_{2}}( y_3^*)^{\beta_{2}}}{2^{\alpha_{2}+\beta_{2}}\alpha_{1}!\alpha_{2}!\beta_{1}!\beta_{2}!} \right]
\end{align*}
Since
$(x_1^*, x_2^*, y_1^*, y_2^*, y_3^*)$ solve the polynomial equations \eqref{eq:symmetric_system_equalities}, 
we have $B_{n,\ell} = 0$ for all $\ell = 1, \dots, \rbar-1$
and
\begin{align}
\label{eq:pf_minimax_asym_maxcoefs}
\max\Big\{|B_{n,\ell}|: \rbar \leq \ell \leq 2\rbar-2 \Big\} \lesssim \epsilon_n^{\rbar} = \frac 1 n.
\end{align}
We therefore have,\begin{align}
\nonumber
h^2\Big(g( \cdot,\bfeta^{(1)}_n), 
 g(\cdot,\bfeta_n^{(2)})\Big) 
 &= \int \left[\frac{\sum \limits_{\ell=\rbar}^{2\rbar-2} B_{n,\ell} \frac{\partial^\ell f}{\partial\theta^l} (x, -\theta_n^{(2)}, v_{1,n}^{(2)}) + R_n(x) }{
 { \sqrt{g(x,\bfeta_n^{(1)})} +
  \sqrt{g(x,\bfeta_n^{(2)})} }}\right]^2 dx \\
  \nonumber
 &\lesssim \int \frac{\sum \limits_{\ell=\rbar}^{2\rbar-2} \left[B_{n,\ell} \frac{\partial^\ell f}{\partial\theta^\ell} (x, -\theta_n^{(2)}, v_{1,n}^{(2)})\right]^2 + R_n^2(x)}{
 {\pi f(x,-\theta_n^{(2)}, v_{1,n}^{(2)}) }} dx \\
 \label{eq:pf_minimax_sym_ending}
 &\lesssim \frac  1 n+ \int \frac{R_n^2(x)}{
    {f(x,-\theta_n^{(2)}, v_{1,n}^{(2)}) }}dx,
\end{align}
where the last inequality follows by integrability of 
$\left[\frac{\partial^\ell f}{\partial\theta^\ell} (\cdot, -\theta_n^{(2)}, v_{1,n}^{(2)})\right]^2/f(\cdot,-\theta_n^{(2)}, v_{1,n}^{(2)})$
for Gaussian densities.
Upon bounding the remainder term in \eqref{eq:pf_minimax_sym_ending} 
in a similar way as Theorem \ref{theorem:convergence_asym} , we arrive at
$$h^2\Big(g(\cdot,\bfeta_n^{(1)}), 
 g(\cdot,\bfeta_n^{(2)})\Big) \lesssim \frac 1 n.$$
The claim follows.  
\hfill $\square$

\subsection{Proof of Propositions \ref{proposition:inequality_necessity} and \ref{proposition:symmetric_system}}
\paragraph{PROOF OF PROPOSITION \ref{proposition:inequality_necessity}}
Set $x_1=x_2 = 0$, and let $y_1 = -y_2 = y_3$ for some arbitrary non-zero real number
$y_3 \in \bbR$. Notice that inequality \eqref{eq:symmetric_system_inequalities} of the symmetric system 
is violated for this setting of variables, and the system of polynomial equations \eqref{eq:symmetric_system_equalities} reduces to
\begin{align}
\label{eq:pf_prop_sym_reduction_b}
 \sum \limits_{\alpha_{2},\beta_{2}} 
 \dfrac{(-1)^{\alpha_2} y_1^{\alpha_2 + \beta_2}}{2^{\alpha_{2}+\beta_{2}} \alpha_{2}! \beta_{2}!}
+ \sum_{\alpha_2} \dfrac{ y_{1}^{\alpha_2}}{2^{\alpha_2}(\alpha_2)!} = 0, \quad \ell =1, \dots, r,
\end{align}
where the first summation in the above display is taken over 
all integers $1 \leq \alpha_2 \leq r$ and $0 \leq \beta_1 \leq r - \alpha_2$
satisfying $2(\alpha_2 + \beta_2) = \ell$, and the second summation is taken
over all integers $1 \leq \alpha_2 \leq r$ such that $2\alpha_2 = \ell$. Clearly, both of these summations are empty
when $\ell$ is odd, hence the equality in equation \eqref{eq:pf_prop_sym_reduction_b} 
holds vacuously for all such $\ell$. It thus remains to prove that the left-hand side of equation \eqref{eq:pf_prop_sym_reduction_b}
vanishes for all even integers $\ell\geq 1$. 

For all even integers $\ell \geq 1$, the left-hand side
of equation \eqref{eq:pf_prop_sym_reduction_b} reduces to
\begin{align*} \frac{y_1^{\ell/2}}{2^{\ell/2}} \sum \limits_{\alpha_{2},\beta_{2}} 
 \dfrac{(-1)^{\alpha_2}}{\alpha_{2}! \beta_{2}!}
+  \dfrac{ y_{1}^{\ell/2}}{2^{\ell/2}(\ell/2)!} 
 &=  \frac{y_1^{\ell/2}}{2^{\ell/2}} \left[\sum \limits_{\alpha_{2}=1}^{\ell/2} 
 \dfrac{(-1)^{\alpha_2}}{\alpha_{2}!(\frac \ell 2-\alpha_2)!}
+  \dfrac{1}{ (\ell/2)!} \right] \\
 &=  \frac{y_1^{\ell/2}}{2^{\ell/2}} \sum \limits_{\alpha_{2}=0}^{\ell/2} 
 \dfrac{(-1)^{\alpha_2}}{\alpha_{2}! (\frac \ell 2-\alpha_2)!}
 \\
 &=  \frac{y_1^{\ell/2}}{2^{\ell/2}(\ell/2)!} \sum \limits_{\alpha_{2}=0}^{\ell/2} 
 {\ell/2 \choose \alpha_2}(-1)^{\alpha_2} \\
 &= 0,
 \end{align*}
by the Binomial Theorem. It follows that $(x_1, x_2, y_1, y_2, y_3)$ solves the system of polynomial
equations in equation  \eqref{eq:symmetric_system_equalities}, thus proving the claim.\hfill $\square$

\paragraph{PROOF OF PROPOSITION \ref{proposition:symmetric_system}}

To prove the claim, we begin by proving that $\ordersym \leq 4$. 
It suffices to show that the system of polynomial equalities
and inequalities \eqref{eq:symmetric_system_equalities} and \eqref{eq:symmetric_system_inequalities}
admits no non-trivial, real-valued, solution when $r = 4$. 
In this case, equalities \eqref{eq:symmetric_system_equalities} read
\begin{align*}
(E_1)    \qquad 0 &= -2x_1 x_2+\frac{y_2}{2} + x_1^2 + \frac{y_1}{2}   \\
(E_2)   \qquad 0 &= x_2y_2-\frac{x_1y_3}{2}-2x_1 x_2^2-\frac{x_1 y_2}{2} +
x_1^2 x_2 +\frac{x_1y_1}{2} \\
(E_3)  \qquad 0 &= \frac{y_2y_3}{4}
+ x_2^2 y_2
+ \frac{y_2^2}{8}
- x_1x_2y_3
-\frac{4 x_1x_2^3}{3}
- x_1x_2y_2
\\ &+ \frac{x_1^2 y_3}{4}
+ x_1^2 x_2^2
+ \frac{x_1^2 y_2}{4}
-\frac{x_1^3 x_2 }{3}
+ \frac{x_1^4}{12}
+ \frac{y_1^2}{8}
+ \frac{x_1^2y_1}{4},
\end{align*}
and inequality \eqref{eq:symmetric_system_inequalities} reads
$$(I) \qquad |x_1|^4 + |y_1|^2 + |y_2|^2 \leq |2x_2-x_1|^4 + |y_3 - y_1|^2 + |y_2+y_3|^2.$$
We first claim that any non-trivial solution $(x_1, x_2, y_1, y_2, y_3)$ to
the above equalities and inequalities must satisfy $x_1 \neq 0$. Indeed, if $x_1 = 0$
by way of a contradiction, then equation $(E_1)$ implies $y_1=-y_2$ 
while equation $(E_2)$ reduces to $0 = x_2y_2$. 
It follows that either $y_2=0$ or $x_2=0$. If $y_2=0$, 
then also $y_1 = 0$, which contradicts the non-triviality of the solution. 
It follows that $x_2 = 0$. Equation $(E_3)$ then reads
$$0 = \frac{y_2y_3}{4}
+ \frac{y_2^2}{8}+ \frac{y_1^2}{8},$$
implying $y_1 = -y_2 = y_3.$
To summarize, if $x_1=0$, the only possible non-trivial
solutions to $(E_1), (E_2), (E_3)$ are of the form $(x_1, x_2, y_1, y_2, y_3)=(0, 0, y_1, -y_1, y_1)$
for $y_1 \in \bbR$. No such solution can satisfy inequality $(I)$.
We thus have a contradiction with the
hypothesis $x_1=0$. 

Since $x_1 \neq 0$, define the variables
\begin{equation}
\label{eq:sym_tildes}
\tilde x_2 = x_2/x_1, \quad
  \tilde y_1 = y_1/x_1^2, \quad
  \tilde y_2 = y_2/x_1^2,\quad
  \tilde y_3 = y_3/x_1^2.
\end{equation}
Equations $(E_1),(E_2),(E_3)$ may then be rewritten as
\begin{align*}
(\bar E_1)    \qquad 0 &= -2 \tilde x_2+\frac{\tilde y_2}{2} + 1 + \frac{\tilde y_1}{2}\\
(\bar E_2)   \qquad 0 &= \tilde x_2 \tilde y_2-\frac{\tilde y_3}{2}-2\tilde x_2^2-\frac{\tilde  y_2}{2} +
\tilde x_2+\frac{\tilde y_1}{2}\\
(\bar E_3)  \qquad 0 &=
\frac{\tilde y_2 \tilde y_3}{4}
+ \tilde x_2^2 \tilde y_2
+ \frac{\tilde y_2^2}{8}
- \tilde x_2\tilde y_3
-\frac{4 \tilde x_2^3}{3}
- \tilde x_2\tilde y_2
\\ &+ \frac{\tilde y_3}{4}
+  \tilde x_2^2
+ \frac{\tilde y_2}{4}
-\frac{\tilde x_2}{3}
+ \frac{1}{12}
+ \frac{\tilde y_1^2}{8}
+ \frac{\tilde y_1}{4}.
\end{align*}
Equation $(\bar E_1)$ implies
\begin{equation}
\label{eq:sym_tilde_x2}
\tilde x_2 = \frac 1 4(2 + \tilde y_1 + \tilde y_2),
\end{equation}
which, combined with $(\bar E_2)$, implies
\begin{equation}
\label{eq:sym_tilde_y3}
\tilde y_3 = \tilde y_1-\tilde y_2 + \frac 1 2 (2+\tilde y_1+\tilde y_2) + \frac 1 2 \tilde y_2(2+\tilde y_1+\tilde y_2) - \frac 1 4 (2+\tilde y_1+\tilde y_2)^2=\frac 1 4 \big[2\tilde y_1 - \tilde y_1^2 + \tilde y_2(\tilde y_2-2)\big].
\end{equation}
With these values of $\tilde x_2, \tilde y_3$, equation $(\bar E_3)$ may be simplified to
$$\tilde y_1 + \tilde y_1^3 + \tilde y_2 + \tilde y_2^3=0.$$
Over $\bbR$, the only solution to this equality is given by $\tilde y_2 = -\tilde y_1$. 
By equations \eqref{eq:sym_tilde_x2} and \eqref{eq:sym_tilde_y3}, 
this leads to $\tilde x_2 = 1/2$ and $\tilde y_3 = \tilde y_1$.
Finally, equation \eqref{eq:sym_tildes} then implies that all non-trivial solutions 
to equations $(E_1), (E_2), (E_3)$ must be of the form
$$x_1 = 2x_2, \quad y_1 = -y_2 = y_3.$$
These values do not satisfy inequality $(I)$. 
We conclude that the system of equalities and inequalities \eqref{eq:symmetric_system_equalities} 
and \eqref{eq:symmetric_system_inequalities} admits no non-trivial, real-valued, solution when $r = 4$,
whence $\ordersym \leq 4$. 

We will now argue that $\ordersym \geq 4$. It suffices to show that 
the system of equations and inequalities \eqref{eq:symmetric_system_equalities} and 
\eqref{eq:symmetric_system_inequalities} admits a solution when $r=3$. In this case, the system reduces
to equations $(E_1), (E_2)$, together with the inequality
$$(I') \qquad |x_1|^3 + |y_1|^{3/2} + |y_2|^{3/2} \leq |2x_2-x_1|^3 + |y_3 - y_1|^{3/2} + |y_2+y_3|^{3/2}.$$
Set $x_1=x_2 \in \bbR$ and $y_3 = 0$ to satisfy inequality $(I')$. 
It can then be seen that $(E_1)$ and $(E_2)$ 
are satisfied whenever $y_2 = 2x_1^2 - y_1$, for any $y_1 \in \bbR$. 
We deduce that $\ordersym \geq 4$, and the claim follows.
\hfill $\square$
\subsection{Proofs of Lemmas}
\label{sec:proofs_lemmas_sym}

\paragraph{PROOF OF LEMMA \ref{lemma:case_b2_sym}}
As in the proof of Theorem \ref{theorem:prelim_sym}, we write $\rbar = \ordersym$ for simplicity.
We prove the Lemma by considering three cases.
\paragraph{Case b.2.1:} $|v_{2,n}^{(2)} - v_{1,n}^{(2)}|^{1/2}/|\theta_{n}^{(1)} - \theta_{n}^{(2)}| \to \infty$ as $n \to \infty$. Invoking the assumption of Case b.2 that $|v_{2,n}^{(2)} - v_{1,n}^{(2)}|/\max \left\{|v_{1,n}^{(1)} - v_{1,n}^{(2)}|, |v_{2,n}^{(1)} - v_{2,n}^{(2)}|\right\} \to \infty$, 
it may be verified that 
\begin{eqnarray}
\dfrac{|B_{n,3}|}{|\theta_{n}^{(1)} - \theta_{n}^{(2)}||v_{2,n}^{(2)} - v_{1,n}^{(2)}|} \to \dfrac{1}{4}. \nonumber
\end{eqnarray}
From the formulation of $\overline{D}_{n}$, it is clear that
\begin{eqnarray}
\dfrac{\overline{D}_{n}}{|\theta_{n}^{(1)} - \theta_{n}^{(2)}||v_{2,n}^{(2)} - v_{1,n}^{(2)}|} \to 0. \nonumber
\end{eqnarray}
Therefore, we obtain that $\max \limits_{1 \leq \ell \leq 2\rbar} |B_{n,\ell}|/\overline{D}_{n} \not \to 0$. Additionally, for each $1 \leq |\alpha| \leq \rbar$, as $n$ is sufficiently large, we have
\begin{eqnarray}
& & \hspace{- 3 em} \dfrac{|\theta_{n}^{(1)} - \theta_{n}^{(2)}|^{\alpha_{1}}|v_{2,n}^{(1)} - v_{2,n}^{(2)}|^{\alpha_{2}}\|\overline{R}_{2,n,\alpha}\|_{\infty}}{\max \limits_{1 \leq \ell \leq 2\rbar} |B_{n,\ell}|} \nonumber \\
& & \leq \dfrac{O\parenth{|\theta_{n}^{(1)} - \theta_{n}^{(2)}|^{\alpha_{1}}|v_{2,n}^{(1)} - v_{2,n}^{(2)}|^{\alpha_{2}}(|\theta_{n}^{(2)}|^{\rbar-|\alpha|+\gamma}+|v_{2,n}^{(2)} - v_{1,n}^{(2)}|^{\rbar-|\alpha|+\gamma}}}{|\theta_{n}^{(1)} - \theta_{n}^{(2)}|v_{2,n}^{(2)} - v_{1,n}^{(2)}|}, \nonumber
\end{eqnarray}
which goes to 0 as $n \to \infty$. 
Hence, we eventually have $\|\overline{R}_n \|_{\infty}/\max \limits_{1 \leq \ell \leq 2\rbar} |B_{n,\ell}| \to 0$. 
\paragraph{Case b.2.2:} $|v_{2,n}^{(2)} - v_{1,n}^{(2)}|^{1/2}/|\theta_{n}^{(1)} - \theta_{n}^{(2)}| \not \to \infty$ as $n \to \infty$. 
Under this assumption, we have 
\begin{align}
\max \left\{|v_{1,n}^{(1)} - v_{1,n}^{(2)}|^{1/2}, |v_{2,n}^{(1)} - v_{2,n}^{(2)}|^{1/2}\right\}/|\theta_{n}^{(1)} - \theta_{n}^{(2)}| \to 0. \nonumber
\end{align} 
If $\theta_{n}^{(1)}/\theta_{n}^{(2)} \not \to -1$, then we obtain that 
\begin{align}
|B_{n,2}|/|\theta_{n}^{(2)} - \theta_{n}^{(1)}|^{2} \not \to 0. \nonumber
\end{align} 
Since $\overline{D}_{n}/|\theta_{n}^{(2)} - \theta_{n}^{(1)}|^{2} \to 0$, the previous result implies that $\max \limits_{1 \leq \ell \leq 2\rbar} |B_{n,\ell}|/\overline{D}_{n} \not \to 0$. Furthermore, for each $1 \leq |\alpha| \leq \rbar$, as $n$ is sufficiently large, we have 
\begin{eqnarray}
& & \hspace{ - 3 em} \dfrac{|\theta_{n}^{(1)} - \theta_{n}^{(2)}|^{\alpha_{1}}|v_{2,n}^{(1)} - v_{2,n}^{(2)}|^{\alpha_{2}}\|\overline{R}_{2,n,\alpha}\|_{\infty}}{\max \limits_{1 \leq \ell \leq 2\rbar} |B_{n,\ell}|} \nonumber \\
& & \leq \dfrac{O\parenth{|\theta_{n}^{(1)} - \theta_{n}^{(2)}|^{\alpha_{1}}|v_{2,n}^{(1)} - v_{2,n}^{(2)}|^{\alpha_{2}}(|\theta_{n}^{(2)}|^{\rbar-|\alpha|+\gamma}+|v_{2,n}^{(2)} - v_{1,n}^{(2)}|^{\rbar-|\alpha|+\gamma}}}{|\theta_{n}^{(1)} - \theta_{n}^{(2)}|^{2}}, \nonumber
\end{eqnarray}
which goes to 0 for all $1 \leq |\alpha| \leq \rbar$. Hence, we have 
$\|\overline{R}_n \|_{\infty}/\max \limits_{1 \leq \ell \leq 2\rbar} |B_{n,\ell}| \to 0$.

As a consequence, we only need to consider the scenario that $\theta_{n}^{(1)}/\theta_{n}^{(2)} \to -1$ as $n \to \infty$. Under that setting, we can verify that if $\abss{B_{n,3}}/\abss{\theta_{n}^{(1)}-\theta_{n}^{(2)}}^3 \to 0$, then we have 
\begin{align}
(v_{2,n}^{(2)} - v_{1,n}^{(2)})/\parenth{\theta_{n}^{(1)} - \theta_{n}^{(2)}}^2 \to 0. \nonumber
\end{align}
However, the above limit leads to
\begin{align}
|B_{n,4}|/|\theta_{n}^{(1)} - \theta_{n}^{(2)}|^{4} \to 5/24. \nonumber
\end{align} 
Therefore, $\max\left\{\abss{B_{n,3}},\abss{B_{n,4}}\right\}/\abss{\theta_{n}^{(1)}- \theta_{n}^{(2)}}^4 \not \to 0$ as $n \to \infty$. As $\overline{D}_{n}/|\theta_{n}^{(1)} - \theta_{n}^{(2)}|^{4} \to 0$, the previous result demonstrates that $\max \limits_{1 \leq \ell \leq 2\rbar} |B_{n,\ell}|/\overline{D}_{n} \not \to 0$. Additionally, we also have 
\begin{align}
\|\overline{R}_{2,n,\alpha}\|_{\infty}/\max \limits_{1 \leq \ell \leq 2\rbar} |B_{n,\ell}| 
   \lesssim \|\overline{R}_{2,n,\alpha}\|_{\infty}/|\theta_{n}^{(1)} - \theta_{n}^{(2)}|^{4} \to 0 \nonumber
\end{align} 
for all $1 \leq |\alpha| \leq \rbar$, which eventually leads to $\|\overline{R}_n\|_{\infty}/\max \limits_{1 \leq \ell \leq 2\rbar} |B_{n,\ell}| \to 0$.
The claim follows.
\hfill $\square$

\paragraph{PROOF OF LEMMA \ref{lemma:case_b3_sym}}
As in the proof of Theorem \ref{theorem:prelim_sym}, we write $\rbar = \ordersym$ for simplicity.

Similarly to the proof of Lemma \ref{lemma:case_a3_asym}, 
we only consider the possibility that 
$$\max \left\{|v_{1,n}^{(1)} - v_{1,n}^{(2)}|, |v_{2,n}^{(1)} - v_{2,n}^{(2)}|, |v_{1,n}^{(2)} - v_{2,n}^{(2)}|\right\}/\biggr\{|\theta_{n}^{(1)} - \theta_{n}^{(2)}||\theta_{n}^{(2)}|\biggr\} \not \to \infty$$
as $n \to \infty$ since the proof argument for other possibilities of this term can be 
carried out in the similar fashion. 
According to the previous assumptions, 
we denote $(v_{1,n}^{(1)}-v_{1,n}^{(2)})/\left\{(\theta_{n}^{(2)} - \theta_{n}^{(1)})\theta_{n}^{(2)}\right\} \to \overline{y}_{1}$ and $(v_{2,n}^{(1)} - v_{2,n}^{(2)})/\left\{(\theta_{n}^{(2)} - \theta_{n}^{(1)})\theta_{n}^{(2)}\right\} \to \overline{y}_{2}$ as $n \to \infty$.  
We will demonstrate that 
\begin{align}
\max \limits_{1 \leq \ell \leq 2\rbar} |B_{n,\ell}|/\left\{|\theta_{n}^{(1)} - \theta_{n}^{(2)}||\theta_{n}^{(2)}|^{3}\right\} \not \to 0. \nonumber
\end{align} 
Assume by the contrary that $\max \limits_{1 \leq \ell \leq 2\rbar} |B_{n,\ell}|/\left\{|\theta_{n}^{(1)} - \theta_{n}^{(2)}||\theta_{n}^{(2)}|^{3}\right\} \to 0$. By dividing both the numerator and denominator of $|B_{n,\ell}|/\left\{|\theta_{n}^{(1)} - \theta_{n}^{(2)}||\theta_{n}^{(2)}|^{3}\right\}$ by $|\theta_{n}^{(1)} - \theta_{n}^{(2)}||\theta_{n}^{(2)}|^{l-1}$ ($2 \leq \ell \leq 4$), as $n \to \infty$, we achieve the following system of polynomial equations
\begin{eqnarray}
\overline{y}_{1}+\overline{y}_{2}= 4, \ \overline{y}_{2} = 2, \ \overline{y}_{2} = 4/3, \nonumber
\end{eqnarray}
which cannot hold. Therefore, $\max \limits_{1 \leq \ell \leq 2\rbar} |B_{n,\ell}|/\left\{|\theta_{n}^{(1)} - \theta_{n}^{(2)}||\theta_{n}^{(2)}|^{3}\right\} \not \to 0$. From the formulation of $D_{n}$, it is clear that $\overline{D}_{n}/\left\{|\theta_{n}^{(1)} - \theta_{n}^{(2)}||\theta_{n}^{(2)}|^{3}\right\} \to 0$. As a consequence, $\max \limits_{1 \leq \ell \leq 2\rbar} |B_{n,\ell}|/\overline{D}_{n} \not \to 0$ as $n \to \infty$. 

Furthermore, we have that
\begin{align}
\dfrac{\|\overline{R}_{2,n,\alpha}\|_{\infty}}{|\theta_{n}^{(1)} - \theta_{n}^{(2)}||\theta_{n}^{(2)}|^{3}}=
 \dfrac{O\parenth{|\theta_{n}^{(2)}|^{\rbar-|\alpha|+\gamma}+|v_{2,n}^{(2)} - v_{1,n}^{(2)}|^{\rbar-|\alpha|+\gamma}}}{|\theta_{n}^{(1)} - \theta_{n}^{(2)}||\theta_{n}^{(2)}|^{3}}, \nonumber
\end{align}
which goes to 0 as $n \to \infty$ for all $1 \leq |\alpha| \leq \rbar$. As a consequence, we have 
\begin{align*}
\|R_n\|_{\infty}/\max \limits_{1 \leq \ell \leq 2\rbar} |B_{n,\ell}| \to 0. 
\end{align*} 
The claim follows. \hfill $\square$

\section{Upper Bounds on the Asymmetric Order}
\label{sec:appendix_rasym}
In this Appendix, we provide upper bounds on the asymmetric order $\orderassym(\pi)$
for certain values of $\pi \in (0,1/2)$. 
We begin with a reduction of the asymmetric 
system of polynomals \eqref{eq:asymmetric_system}. 

We claim that the system does not admit a non-trivial solution $(x_1, x_2, y_1, y_2, y_3) \in \bbR^5$ with $x_1=0$
when $r=6$. Indeed, when $x_1 = 0$, the system reduces to 
\begin{align}
\label{eq:system_rasym_ub}
(1-\pi) \sum \limits_{\alpha_{2},\beta_{1},\beta_{2}} 
\dfrac{1}{2^{\alpha_{2}+\beta_{2}}} \dfrac{(c+1)^{\beta_{1}}x_{2}^{\beta_{1}}y_{2}^{\alpha_{2}}y_{3}^{\beta_{2}}}{\alpha_{2}!\beta_{1}!\beta_{2}!} 
+ \pi \sum \limits_{\alpha_{2}} \dfrac{1}{2^{\alpha_{2}}}\dfrac{y_{1}^{\alpha_{2}}}{\alpha_{2}!}  = 0, ~~
 \quad \ell=1,\ldots,r 
\end{align}
where
the first sum is over all nonnegative integers $\alpha_{2},\beta_{1},\beta_{2}$ such that
$\beta_{1}+2\alpha_{2}+2\beta_{2} = \ell$, $1 \leq \alpha_{2} \leq r$, 
and $0 \leq \beta_{1}+\beta_{2} \leq r - \alpha_{2}$,
while, in the second sum, $1 \leq \alpha_{2} \leq r$
ranges over all integers satisfying $\alpha_{2}=\ell/2$.
In particular, the second sum is empty whenever $\ell$ is odd. 

The equation for $\ell=1$ of the system holds trivially. The equation for $\ell=2$
implies $y_2=-cy_1$, while that of $\ell=3$ implies $y_2x_2=0$. 
If $y_2 = 0$, then also $y_1=0$ and the solution becomes trivial, thus it follows
that $x_2 = 0$. The system \eqref{eq:system_rasym_ub} then
reduces to 
\begin{align*}
(1-\pi) \sum_{\alpha_{2},\beta_{2}} 
\dfrac{(-cy_1)^{\alpha_{2}}y_{3}^{\beta_{2}}}{2^{\alpha_{2}+\beta_{2}}\alpha_{2}!\beta_{2}!} 
+ \pi \sum \limits_{\alpha_{2}} \dfrac{y_{1}^{\alpha_{2}}}{2^{\alpha_2}\alpha_{2}!}=0, \quad \ell=1, \dots, r.
\end{align*}
By definition of the ranges in the above summations, both summations are empty when $\ell$ is an odd integer.
When $\ell$ is even, the above display reduces to
\begin{align*}
(1-\pi) \sum_{\alpha_{2}=1}^{\ell/2} 
\dfrac{(-cy_1)^{\alpha_{2}}y_{3}^{\ell/2-\alpha_2}}{2^{\ell/2}\alpha_{2}!(\ell/2-\alpha_2)!} 
+ \pi \dfrac{y_{1}^{\ell/2}}{2^{\ell/2}(\ell/2)!}=0, \quad \ell/2 = 1, \dots, \lfloor r/2\rfloor.
\end{align*}
Taking $\ell=4$  implies 
$y_3 = -\frac{(1+c)y_1}{2c}$, thus we have the further reduction
\begin{align*}
(1-\pi) \sum_{\alpha_{2}=1}^{\ell/2} 
\dfrac{(-c)^{\alpha_{2}}(-\frac{(1+c)}{2c})^{\ell/2-\alpha_2}y_1^{\ell/2}}{2^{\ell/2}\alpha_{2}!(\ell/2-\alpha_2)!} 
+ \pi \dfrac{y_{1}^{\ell/2}}{2^{\ell/2}(\ell/2)!}=0, \quad \ell/2 = 1, \dots, \lfloor r/2\rfloor.
\end{align*}
Since $y_1 \neq 0$, the above display no longer depends on the variables $x_1, x_2, y_1, y_2, y_3$, and reduces to
\begin{align*}
(1-\pi) \sum_{\alpha_{2}=1}^{\ell/2} 
\dfrac{(-\frac{(1+c)}{2})^{\ell/2-\alpha_2}c^{\ell/2}}{\alpha_{2}!(\ell/2-\alpha_2)!} 
+ \pi \dfrac{1}{(\ell/2)!}=0, \quad \ell/2 = 1, \dots, \lfloor r/2\rfloor.
\end{align*}
It can be seen by direct verification that the above display does not hold when $\ell = 6$
provided  $\pi\neq 1/2$. 
Therefore, there exists no non-trivial solution to the asymmetric system with $r=6$ when $x_1 = 0$. 
In what follows, we will therefore assume $x_1 \neq 0$, and show that the system continues to have no solution
for $r=6$ for a range of values of $\pi$.

Fix $r=6$. Since $x_1 \neq 0$, each equation of the the asymmetric system \eqref{eq:asymmetric_system} may be divided by $x_1^\ell$,
leading to the system
\begin{align}
\label{eq:final_sys_grobner}
(1-\pi) \sum \limits_{\alpha_{1},\alpha_{2},\beta_{1},\beta_{2}} \dfrac{c^{\alpha_{1}}(c+1)^{\beta_{1}}(-1)^{\alpha_{1}}\tilde x_{2}^{\beta_{1}}\tilde y_{2}^{\alpha_{2}}\tilde y_{3}^{\beta_{2}}}{2^{\alpha_{2}+\beta_{2}}\alpha_{1}!\alpha_{2}!\beta_{1}!\beta_{2}!} 
+ \pi \sum \limits_{\alpha_{1},\alpha_{2}} \dfrac{1}{2^{\alpha_{2}}}\dfrac{\tilde 
y_{1}^{\alpha_{2}}}{\alpha_{1}!\alpha_{2}!}  = 0, \quad \ell=1,\ldots,6,
\end{align}
where $\tilde x_2 = x_2/x_1$ and $\tilde y_j = y_j/x_1^2$ for $j=1,2,3$. 
We compute\footnote{See \url{https://github.com/tmanole/Gaussian-mixture-twocomp}.} a
reduced Gr\"obner basis of the above polynomials over $\bbC$,
in the Mathematica programming language \citep{wolfram1999},
for $\pi \in \{i/100: 1 \leq i \leq 49, i \in \bbN\}$. For all such values of $\pi$,
we obtain the Gr\"obner basis $\{1\}$. It follows that the system of equations \eqref{eq:final_sys_grobner}
does not admit any solution for these values of $\pi$. Together with the result of Proposition \ref{proposition:asymmetric_system},
we conclude $\orderassym(\pi) = 6$ for all $\pi \in  \{i/100: 1 \leq i \leq 49, i \in \bbN\}$.
 
\section{Numerical Supplement}
\label{sec:appendix_numerical}

\subsection{Simulation Specifications}
In this Appendix, we provide additional details for the numerical experiments
in Section \ref{sec:simulations}. 

The specific form of the EM algorithm for model \eqref{eqn:general_model} is straightforward
to derive, and is summarized in Algorithm \ref{alg:EM}. In our experiments, we use the convergence
criterion $\epsilon=10^{-8}$, and we halted the EM algorithm if its number of iterations
exceeded $T = 2,000$. 

Since the purpose of our simulations is to illustrate the theoretical rate of convergence
of the parameters in location-scale Gaussian mixtures, we initialize the EM algorithm based
on the true parameter values. Specifically, we initialize
the location and scale parameters respectively 
by uniformly sampling from the intervals $[\theta_n -  n^{-1/14}, \theta_n + n^{-1/14}]$, 
and  $[v_{j,n} - n^{-1/7}, v_{j,n} + n^{-1/7}]$ for $j=1, 2$.
Here $\theta_n, v_{1,n}, v_{2,n}$ denote the true parameters under each of Models 1 and 2.
For each replication in our simulations, 
we run the EM algorithm five times with distinct starting values of this form, and retain the fitted solution
which achieved the highest likelihood.

\begin{algorithm}[t]
   \caption{\label{alg:EM} EM Algorithm for Model \eqref{eqn:general_model}}
\begin{algorithmic}
   \STATE {\bfseries Input:} Sample $Y_1, \dots, Y_n$; 
                             Starting values  $\bfeta_n^{(0)} = (\theta_n^{(0)}, v_{1,n}^{(0)},v_{2,n}^{(0)})$;
                             Mixing proportion $\pi \in (0,1/2]$  			
   \STATE {\bfseries Output:} Approximate maximum likelihood estimate of $\bfeta_n$
   \STATE $t \leftarrow 0; \ c \leftarrow \pi/(1-\pi)$
   \WHILE {$|\ell_n(\bfeta_n^{(t)}) - \ell_n(\bfeta_n^{(t+1)})| > \epsilon$}
   \STATE 1. Let
   $$w_{i}^{(t)} \leftarrow \frac{\pi f(Y_i; -\theta_n^{(t)}, v_{1,n}^{(t)})}
   						    {\pi f(Y_i; -\theta_n^{(t)}, v_{1,n}^{(t)})+ 
    						 (1-\pi) f(Y_i; c\theta_n^{(t)}, v_{2,n}^{(t)})}, \quad i=1, \dots, n.$$
   \STATE 2. Update $\theta_n^{(t)}$,
   $$\theta_n^{(t+1)} \leftarrow \frac 1 {cn + (1-c)\sum_{i=1}^n w_i} 
   								  \sum_{i=1}^n \Big[c(1-w_i^{(t)})Y_i - w_i^{(t)}Y_i\Big]$$
   \STATE 3. Update $v_{1,n}^{(t)}, $
   $$v_{1,n}^{(t+1)} \leftarrow \frac 1 {\sum_{i=1}^n w_i^{(t)}} \sum_{i=1}^n w_i\Big(Y_i + \theta_n^{(t+1)}\Big)^2.$$
   \STATE 4. Update $v_{2,n}^{(t)}$,
   $$v_{2,n}^{(t+1)} \leftarrow \frac 1 {  n- \sum_{i=1}^n w_i^{(t)}} \sum_{i=1}^n  \Big(1-w_i^{(t)}\Big)\Big(Y_i-c\theta_n\Big)^2.$$ 								  				 
   \STATE 1. Update $Y_{j}^{(t)}$ and $b_{j}^{(t)}$ for $1 \leq j \leq m$:
   \IF{$t \geq T$}
   \STATE \textbf{break}
   \ENDIF
	$t \leftarrow t+1$
   \ENDWHILE
\end{algorithmic}
\end{algorithm}

\clearpage
\subsection{Additional Simulation Results}
\begin{figure}[hbp!]
\centering
\begin{subfigure}{0.49\textwidth}
  \includegraphics[width=.92\linewidth]{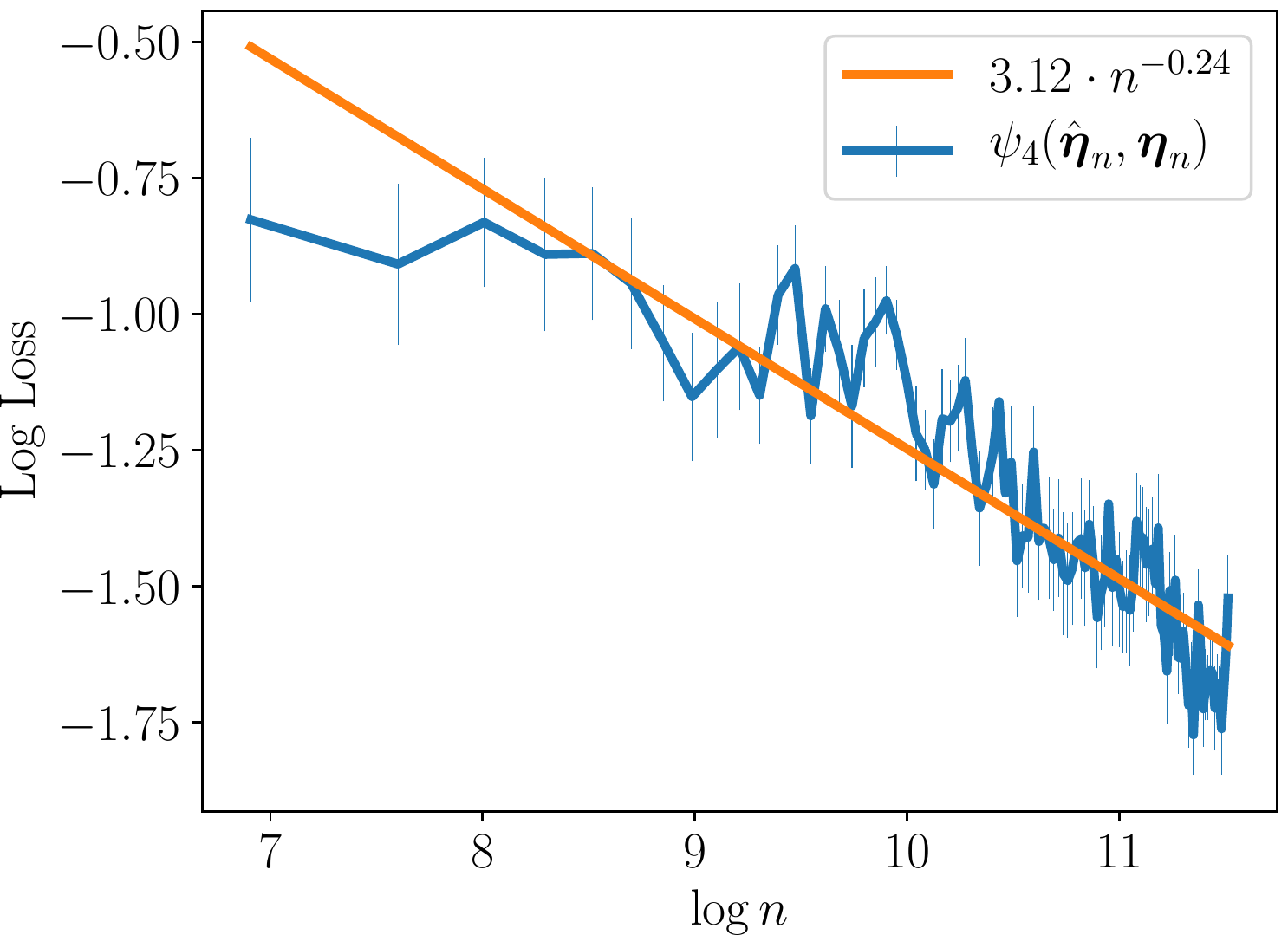}
  \caption{Model A, $\pi=0.5$.}
\end{subfigure}%
\begin{subfigure}{0.49\textwidth}
  \centering
  \includegraphics[width=.92\linewidth]{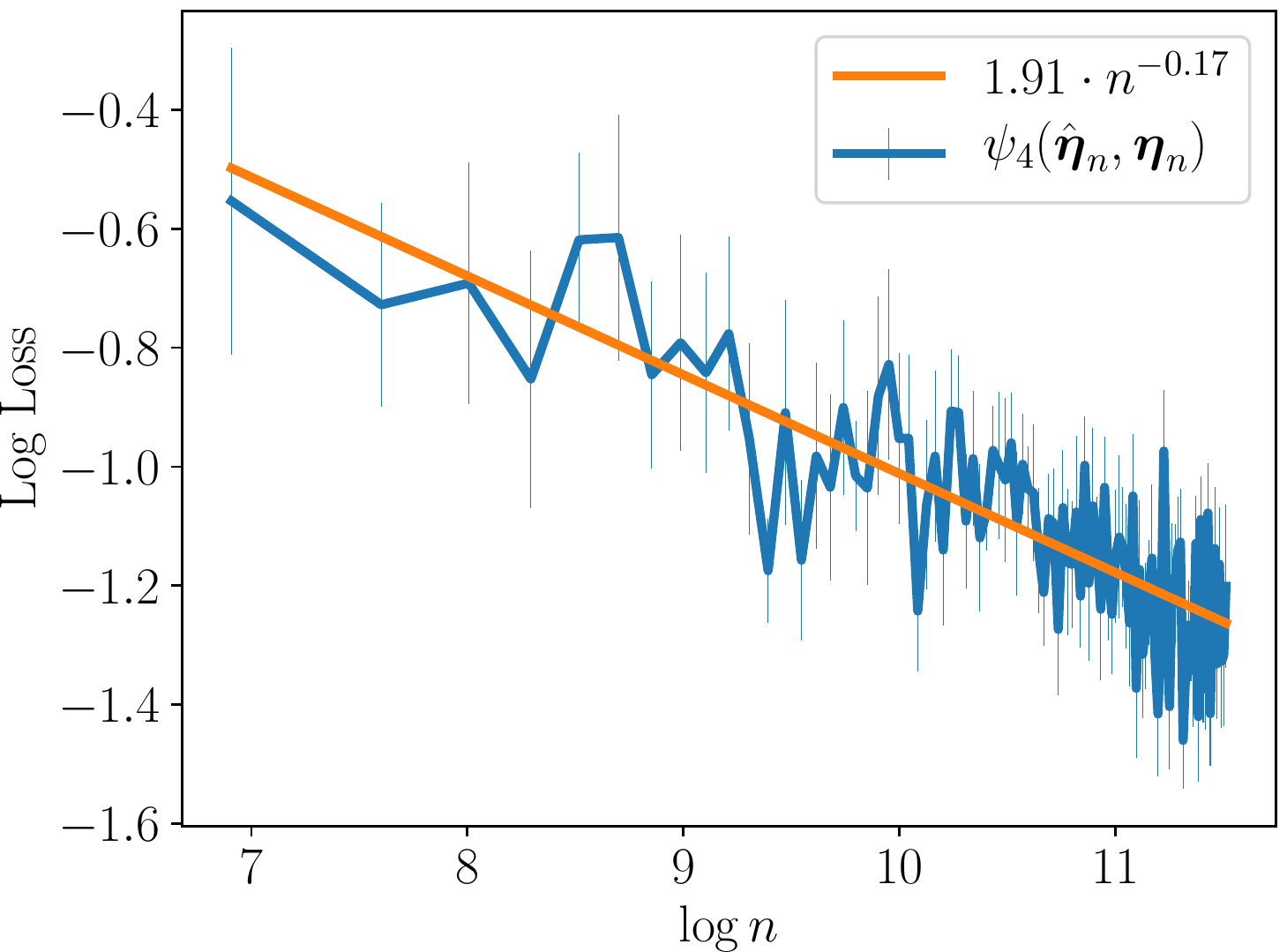}
  \caption{Model S', $\pi=0.5$.}
\end{subfigure}%
\caption{\label{fig:appendix} Log-log scale plots for the simulation results under Model A with $\pi=1/2$ and under Model S'.
For each model and sample size, the MLE $\hbfeta_n$ is computed on 10 independent samples, and its average distance from
$\bfeta_n$ is plotted with error bars representing one empirical standard deviation. 
}
\end{figure}

We now provide two additional simulation results in the symmetric regime.
In Figure \ref{fig:appendix}(a), we report the result of Model A from Section 
\ref{sec:simulations} under $\pi=1/2$. Furthermore, in Figure \ref{fig:appendix}(b), 
we report the results under the following distinct parameter setting
$$\textbf{Model S'}: \quad \theta_n = n^{-1/8}, \quad
  v_{1,n} = 1 + n^{-1/4}/3, \quad
  v_{2,n} = 1 + n^{-1/4}  /6.$$
It can be seen that these parameter settings respectively achieve the approximate $n^{-1/4}$ 
and $n^{-1/6}$ empirical rates of convergence.

\clearpage

\end{document}